\def\year{2022}\relax
\documentclass[letterpaper]{article} 

\usepackage{hyperref}

\usepackage{graphicx} 
    
\usepackage{aaai22}  
\usepackage{times}  
\usepackage{helvet}  
\usepackage{courier}  

\usepackage{booktabs,amsfonts,microtype,enumitem,bm,pifont,amsmath,amssymb,amsthm,cleveref, multicol}

\usepackage{natbib}  
\usepackage{caption} 
\DeclareCaptionStyle{ruled}{labelfont=normalfont,labelsep=colon,strut=off} 
\usepackage{algorithm}
\usepackage{algorithmic}

\usepackage{newfloat}
\usepackage{listings}
\floatstyle{ruled}
\newfloat{listing}{tb}{lst}{}
\floatname{listing}{Listing}
\newcommand  \tmop[1]    {{\ensuremath{\operatorname{#1}}}}

\def \a {{\boldsymbol{a}}}
\def \b {{\boldsymbol{b}}}

\def \e {{\boldsymbol{e}}}
\def \f {{\boldsymbol{f}}}
\def \d {{\boldsymbol{d}}}
\def \g {{\boldsymbol{g}}}
\def \h {{\boldsymbol{h}}}

\def \u {{\boldsymbol{u}}}
\def \v {{\boldsymbol{v}}}
\def \w {{\boldsymbol{w}}}
\def \x {{\boldsymbol{x}}}
\def \y {{\boldsymbol{y}}}

\def \bxi    {{\boldsymbol {\xi}}}
\def \vlambda {{\boldsymbol{\lambda}}}

\def \eye  {{\boldsymbol{I}}}
\def \minimize  {\operatorname*{minimize}}

\def \st        {\operatorname*{subject\ to\ }}
\def \nn        {\nonumber}
\def \R {\mathbb{R}}

\def \S {\mathbb{S}}

\def \eye       {\mathbf{I}}
\def \zero      {\mathbf{0}}
\def \one       {\mathbf{1}}
\def \lg        {\langle}
\def \rg        {\rangle}
\def \diag      {\tmop{diag}}

\def \rank      {\tmop{rank}}

\def \dist      {\tmop{dist}}

\def \polylog {\tmop{polylog}}

\newtheorem{lem}{Lemma}
\newtheorem{theorem}{Theorem}
\newtheorem{cor}{Corollary}
\newtheorem{definition}{Definition}
\newtheorem{assumption}{Assumption}
\theoremstyle{remark}

\title{Local and Global Convergence of General Burer-Monteiro Tensor Optimizations}
\author{
   Shuang Li,\textsuperscript{\rm 1}
   Qiuwei Li \textsuperscript{\rm 2}
   }

\affiliations {
     \textsuperscript{\rm 1} Department of Mathematics, University of California, Los Angeles \\
     \textsuperscript{\rm 2} Alibaba Group (US), Damo Academy 
     \\
     shuangli@math.ucla.edu, liqiuweiss@gmail.com.com
}

\usepackage{bibentry}

\begin{document}

\maketitle

\begin{abstract}
Tensor optimization is crucial to massive machine learning and signal processing tasks. In this paper, we consider tensor optimization with a convex and well-conditioned objective function and reformulate it into a nonconvex optimization using the Burer-Monteiro type parameterization. We analyze the local convergence of applying vanilla gradient descent to the factored formulation and establish a local regularity condition under mild assumptions. We also provide a linear convergence analysis of the gradient descent algorithm started in a neighborhood of the true tensor factors.  Complementary to the local analysis, this work also characterizes the global geometry of the best rank-one tensor approximation problem and demonstrates that for orthogonally decomposable tensors the problem has no spurious local minima and all saddle points are strict except for the one at zero which is a third-order saddle point.
\end{abstract}

\section{Introduction}
\label{sec:intro}
Tensors, a multi-dimensional generalization of vectors and matrices, provide natural representations for multi-way datasets and find numerous applications in machine learning and signal processing, including video processing~\cite{liu2012tensor}, hyperspectral imaging~\cite{li2015low,sun2020weighted}, collaborative filtering \cite{overcomplete:Chen:2005jn}, latent graphical model learning \cite{overcomplete:anandkumar2017analyzing}, independent component analysis (ICA) \cite{overcomplete:Cardoso:1989eo},   dictionary learning~\cite{barak2015dictionary},  neural networks compression ~\cite{phan2020stable,bai2021decomvqanet},  Gaussian mixture estimation \cite{overcomplete:Hsu:2013iz},  and psychometrics \cite{overcomplete:Smilde:2005wf}. See \cite{sidiropoulos2017tensor} for a review. All these applications involve solving certain optimizations over the space of low-rank tensors:
\begin{equation}
\begin{aligned}
\minimize_{T} &~ f(T) ~~
\st~\rank(T)\le r.
\label{eqn:f(T)}
\end{aligned}
\end{equation}
Here $f(\cdot)$ is a problem dependent objective function with tensor argument and $\rank(\cdot)$ calculates the tensor rank.  The rank of matrices is well-understood and has many equivalent definitions, such as the dimension of the range space, or the size of largest non-vanishing minor, or the number of nonzero singular values. The latter is also equal to the smallest number of rank-one factors that the matrix can be written as a sum of. The tensor rank, however, has several non-equivalent variants, among which the Tucker rank \cite{overcomplete:Kolda:2009dh} and the Canonical Polyadic (CP) rank \cite{grasedyck2013literature} are most well-known. The CP tensor rank is a more direct generalization from the matrix case and is precisely equal to the minimal number of terms in a rank-one tensor decomposition. It is also the preferred notion of rank in applications. Unfortunately, while the Tucker rank can be found by performing the  higher-order singular value decomposition (HOSVD) of the tensor, the CP rank is NP-hard to compute \cite{overcomplete:Hillar:2013by}. Even though some recent works \cite{yuan2016tensor,barak2016noisy, li2016super, li2017convex, li2015overcomplete, tang2015guaranteed} study the convex relaxation methods based on the tensor nuclear norm, which is also NP-hard to compute \cite{overcomplete:Hillar:2013by}. Therefore, this work seeks alternative ways to solve the CP rank-constrained tensor optimizations.

 \subsection{General Burer-Monteiro Tensor Optimizations}
Throughout this paper, we focus on third-order, {\em symmetric} tensors and assume that $f:\R^{n\times n\times n}\rightarrow \R$ is a general convex function and has a unique global minimizer $T^\star$ that admits the following (symmetric-)rank-revealing decomposition:   
\begin{align}
T^{\star}=\sum_{p=1}^{r} {c_p^\star} \widehat{\u}_{p}\otimes\widehat{\u}_{p} \otimes \widehat{\u}_{p} \in \R^{n \times n \times n},
\label{eqn:true:tensor0}
\end{align}
where  $\widehat{\u}_{p}$'s are the normalized tensor factors living on the unit spheres $\mathbb{S}^{n-1}$ and $c_p^\star$'s are the decomposition coefficients. Without loss of generality, we can always assume $c_p^\star>0$, since otherwise we can absorb its sign into the normalized tensor factors.

Note that the global optimal tensor in \eqref{eqn:true:tensor0} can be rewritten as
\begin{align}
T^{\star}
&=\sum_{p=1}^{r}  ({c_p^\star}^{1/3}\widehat{\u}_{p}) \otimes({c_p^\star}^{1/3}\widehat{\u}_{p}) \otimes ({c_p^\star}^{1/3}\widehat{\u}_{p})
\nn
\\
&\doteq U^\star\circ U^\star \circ U^\star,
\label{eqn:truemodel}	
\end{align}
where $U^\star\doteq\big[{c_1^\star}^{1/3}\widehat{\u}_{1}~{c_2^\star}^{1/3}\widehat{\u}_{2}~\cdots~{c_r^\star}^{1/3}\widehat{\u}_{r}\big]$ can be viewed as the ``cubic root'' of $T^\star$. Noting that the ``cubic-root'' representation \eqref{eqn:truemodel} has permutation ambiguities, that is, different columnwise permutations of $U^\star$ would generate the same tensor in \eqref{eqn:truemodel}:
$[\u_{i_1}^\star~\u_{i_2}^\star~\cdots~\u_{i_r}^\star]\circ [\u_{i_1}^\star~\u_{i_2}^\star~\cdots~\u_{i_r}^\star]
\circ [\u_{i_1}^\star~\u_{i_2}^\star~\cdots~\u_{i_r}^\star]
=U^\star\circ U^\star\circ U^\star$	
for any permutation $(i_1,i_2,\cdots,i_r)$ of the index $(1,2,\cdots,r).$
This immediately implies that $U^\star$ and its columnwise permutations all give rise to global minimizers of the following reformulation of the optimization \eqref{eqn:f(T)}:
\begin{align}
\minimize_{U \in \R^{n \times r}} ~ f(U\circ U\circ U). 
\label{eqn:paramter}
\end{align}
Note that this new factorized formulation has explicitly encoded the rank constraint $\rank(T)\le r$ into the factorization representation $T = U\circ U\circ U$. As a result, the rank-constrained optimization problem \eqref{eqn:f(T)} on tensor variables reduces to the above unconstrained optimization of matrix variables, avoiding dealing with the difficult rank constraint at the price of working with a highly non-convex objective function in $U$. Indeed, while the resulting optimization \eqref{eqn:paramter} has no rank constraint, a smaller memory footprint, and is more amenable for applying simple iterative algorithms like gradient descent, the permutational invariance of $f(U\circ U\circ U)$ implies that saddle points abound the optimization landscape among the exponentially many equivalent global minimizers. Unlike the original convex objective $f(T)$ that has an algorithm-friendly landscape where all the stationary points correspond to the global minimizers, the landscape for the resulting nonconvex formulation $f(U\circ U\circ U)$ is not well-understood. On the other hand, simple local search algorithms applied to \eqref{eqn:paramter} has exhibited superb empirical performance. As a first step towards understanding of the power of using the factorization method to solve tensor inverse problems, this work will focus on  characterizing the local convergence of applying vanilla gradient descent to  the general problem \eqref{eqn:paramter}, as well as the global convergence of a simple variant.

\subsection{Related Work}
\subsubsection{Burer-Monteiro  Parameterization Method}
The idea of transforming the rank-constrained problem into an unconstrained problem using explicit factorization like $T=U\circ U\circ U$ is pioneered by Burer and Monteiro \cite{burer2003nonlinear,burer2005local} in solving matrix optimization problems with a rank constraint 
\begin{equation} 
\begin{aligned}	
\minimize_{X\in\R^{n\times n}} &~~f(X)
\\
\st&~~ \rank(X)\le r~\mathrm{and}~X\succeq 0
\end{aligned}
\end{equation}
To deal with the rank constraint as well as the positive semidefinite constraint, the authors there proposed to firstly factorize a low-rank matrix $X = UU^\top$  with $U \in \R^{n\times r}$ and $r$ chosen according to the rank constraint. Consequently, instead of minimizing an objective function $f(X)$ over all symmetric, positive semidefinite matrices of rank at most $r$, one can focus on an unconstrained nonconvex optimization: 
\[
\minimize_{U\in\R^{n\times r}} f(UU^\top).
\]
Inspired by \cite{burer2003nonlinear,burer2005local},  an intensive research effort has been devoted to investigating the theoretical properties of this factorization/parametrization method \cite{ge2016matrix, ge2017no, park2017non,chi2019nonconvex,li2018non, zhu2018global, zhu2021global, li2017geometry, zhu2019distributed, li2020global}. In particular, by analyzing the landscape of the resulting optimization, many authors have found that various low-rank matrix recovery problems in factored form--despite nonconvexity--enjoy a favorable landscape where all second-order stationary points are global minima.

%

\subsubsection{Tensor Decomposition and Completion} Another line of related work is nonconvex tensor factorization/completion. When the convex objective function $f(T)$ in \eqref{eqn:f(T)} is the squared Euclidean distance between the tensor variable $T$ and the ground-truth tensor $T^\star$, i.e., $f(T)=\|T-T^\star\|_F^2$, the resulting factorized problem \eqref{eqn:paramter} reduces to a (symmetric) tensor decomposition problem:
\begin{align}
\minimize_{U \in \R^{n \times r}} ~ f(U\circ U\circ U)=\|U\circ U\circ U-T^\star\|_F^2. 
\label{eqn:tendec}
\end{align}
Tensor decomposition aims to identify the unknown rank-one factors from available tensor data. This problem is the backbone of several tensor-based machine learning methods, such as independent component analysis \cite{overcomplete:Cardoso:1989eo} and   collaborative filtering \cite{overcomplete:Chen:2005jn}. Unlike the similarly defined matrix decomposition, which has a closed-form solution given by the singular value decomposition, 
the tensor decomposition solution generally has no analytic expressions  and is NP-hard to compute in the worst case~\cite{overcomplete:Hillar:2013by}.  When the true tensor $T^\star$ is a fourth-order symmetric orthogonal tensor, i.e., there is an orthogonal matrix $U^\star$ such that $T^\star=U^\star \circ U^\star \circ U^\star \circ U^\star$, Ge et al. \cite{ge2015escaping} designed a new objective function
\begin{equation*}
\begin{aligned}
\widetilde f (U)\doteq \sum_{i\neq j} \lg T^\star, \u_i\otimes \u_i \otimes \u_j \otimes \u_j\rg
\end{aligned}
\end{equation*}
 and showed that, despite its non-convexity,  the objective function $\widetilde f (U)$ has a benign landscape on the sphere where all the local minima are global minima and all the saddle points have a Hessian with at least one negative eigenvalue. Later, \cite{qu2019geometric} relax the orthogonal condition to near-orthogonal condition, resulting to landscape analysis to fourth-order overcomplete tensor decomposition. The work \cite{ge2015escaping} has spurred many followups that dedicate on the analysis of the nonconvex optimization landscape of many other problems \cite{ge2016matrix, ge2017no,bhojanapalli2016global, park2017non,chi2019nonconvex}. The techniques developed in \cite{ge2015escaping}, however, are not directly applicable to solve the original rank-constrained tensor optimization problem \eqref{eqn:tendec}. In addition, \cite{ge2015escaping} mainly considered fourth-order tensor decomposition, which cannot be trivially extended to analyze other odd-order tensor decompositions. More recently, Ge and Ma \cite{ge2017optimization} studied the problem of maximizing  
 \begin{equation*}
\begin{aligned}
\widehat f(\u)= \lg  T, \u\otimes\u\otimes\u\otimes \u\rg	
\end{aligned}
\end{equation*}
 on the unit sphere and presented a local convergence of applying vanilla gradient descent to  the problem. Although this formulation together with iterative rank-1 updates lead to algorithms with convergence guarantees for tensor decomposition, it is not flexible enough to deal with general rank-constrained problem \eqref{eqn:f(T)}. Similar rank-1 updating methods for tensor decomposition have also been investigated in~\cite{overcomplete:anandkumar2017analyzing,anandkumar2015learning, anandkumar2014guaranteed,anandkumar2014tensor}.

More recently, \cite{chi2019nonconvex,cai2021nonconvex2} apply the factorization formulation to the tensor completion problem and focuses on solving 
\begin{align}
\minimize_{U \in \R^{n \times r}} ~ \left\|P_{\Omega}\left(U\circ U\circ U-T^\star\right)\right\|_F^2,
\label{eqn:tencp}
\end{align}
where  $P_{\Omega}$ is the the orthogonal projection of any tensor $T$ onto the subspace indexed by the observation set $\Omega$. 
\cite{chi2019nonconvex,cai2021nonconvex2} proposed a vanilla gradient descent following a rough initialization and proved the vanilla gradient descent  could faithfully complete the tensor and retrieve all individual tensor factors within nearly linear time when the rank $r$ does not exceed $O(n^{1/6})$. Compared with these prior state of the arts, our convergence analysis improves the order of rank $r$ and extends the focus to general cost functions.

\subsection{Main Contributions and Organization}
To solve the rank-constrained tensor optimization problem \eqref{eqn:f(T)}, we directly work with the  Burer-Monteiro factorized formulation  \eqref{eqn:paramter} with a general convex function $f(\cdot)$ and focus on solving ~\eqref{eqn:paramter} using (vanilla) gradient descent
\begin{align}
U^+=U-\eta \nabla_U f(U\circ U\circ U),
\label{eqn:algorithm}
\end{align}
where $U^+$ is the updated version of the current variable $U$, $\eta$ is the stepsize that will be carefully tuned to prevent gradient descent  from diverging, and
$\nabla_U f$ is the gradient of $f(U\circ U\circ U)$ with respect to $U$.

In this work, we show that the factorized tensor minimization problem~\eqref{eqn:paramter} satisfies the local regularity condition under certain mild assumptions. With this local regularity condition, we further prove a linear convergence of the gradient descent algorithm in a neighborhood of true tensor factors. In particular, we have shown that solving the factored tensor minimization problem~\eqref{eqn:paramter} with gradient descent~\eqref{eqn:algorithm} is guaranteed to identify the target tensor $T^\star$ with high probability if $r = O\big(n^{1.25}\big) $ and $n$ is sufficiently large. This implies that we can even deal with the scenario where the rank of the target tensor $T^\star$ is larger than the individual tensor dimensions, the so called overcomplete regime that are considered challenging to tackle in practice.

Finally, as a complement to the local analysis, we study the global landscape of best rank-1 approximation of a third-order orthogonal tensor and we show that this problem has no spurious local minima and all saddle points are strict saddle points except for the one at zero, which is a third-order saddle point.
 

\subsubsection{Organization}The remainder of this work is organized as follows. In Section~\ref{sec:main}, we first briefly introduce some basic definitions and concepts used in tensor analysis and then present the local convergence of applying vanilla gradient descent to the tensor minimization problem~\eqref{eqn:paramter} and provide a linear convergence analysis for the gradient descent algorithm~\eqref{eqn:algorithm}. In Section~\ref{sec:global}, we switch to analyze the global landscape of orthogonal tensor decomposition. Numerical simulations are conducted in Section~\ref{sec:simu} to further support our theory. Finally, we conclude our work in Section~\ref{sec:conc}.



\section{Local Convergence}

\label{sec:main}


In this section, we first briefly review some fundamental concepts and definitions in tensor analysis. A tensor with order higher than $3$ can be viewed as a high-dimensional extension of vectors and matrices. In this work, we mainly focus on the third-order symmetric tensors. Any such tensor admits symmetric rank-one decompositions of the following form:  
\begin{align*}
T=\sum_{p=1}^{r} {c_p} \u_{p}\otimes\u_{p} \otimes \u_{p} \in \R^{n \times n \times n}
\end{align*}
with $\|\u_{p}\|_2 = 1$ and $c_p > 0$, $1\leq p \leq r$. The above decomposition is also called the {\em Canonical Polyadic} (CP) decomposition of the tensor $T$ \cite{hong2020generalized}. The minimal number of factors $r$ is defined as the {\em (symmetric) rank} of the tensor $T$. Denote $T(i_1,i_2,i_3)$ as the $(i_1,i_2,i_3)$-th entry of a tensor $T$. We define the {\em inner product} of any two tensors $X,Y\in\R^{n \times n\times n}$ as $\langle X, Y \rangle \doteq 	\sum_{i_1,i_2,i_3 = 1}^n X(i_1,i_2,i_3) Y(i_1,i_2,i_3)$. The induced {\em Frobenius norm} of a tensor $T$ is then defined as
$\|T\|_F \doteq \sqrt{\langle T,T \rangle}.$ For a tensor $T\in\R^{n \times n\times n}$, we denote its unfolding/matricization along the first dimension as $T_{(1)}=[T(:,1,1)~T(:,2,1)~\cdots~T(:,n,n)]\in\R^{n \times n^2}.$

\bigskip

We proceed to present the local convergence of applying vanilla gradient descent to  the factored tensor minimization problem~\eqref{eqn:paramter}. Before that, we introduce several  definitions used throughout the work.

\begin{definition}
A function $f: \R^{n \times n\times n} \rightarrow \R$ is $(r,m,M)$-restricted strongly convex and smooth if 
\begin{align*}
m\|Y-X\|_F\leq \|\nabla f(Y)-\nabla f(X)\|_F \leq M\|Y-X\|_F
\end{align*}
holds for any symmetric tensors $X,Y\in\R^{n \times n\times n}$ of rank at most $r$ with some positive constants $m$ and $M$.
\end{definition}
For example, $f(T)=\frac{1}{2}\|T-T^\star\|_F^2$ is such a $(r,m,M)$-restricted strongly convex and smooth function for arbitrary $r\in \mathbb{N}$ with $M=m=1$, and its global minimizer is $T = T^\star$.

\begin{definition} The distance between two factored matrices $U_1$ and $U_2$ is defined as
\begin{align*}
\dist(U_1,U_2)&=\min_{\tmop{Permutation} P}\|U_1-U_2P\|_F.
\end{align*}
\end{definition}
Denote 
\begin{equation}
P_{U_1}=\arg\min_{\tmop{Permutation} P}\|U_1-U_2P\|_F. 
\end{equation}
Then, we can rewrite the distance between $U_1$ and $U_2$ as
\begin{equation}
\dist(U_1,U_2)=\|U_1-U_2P_{U_1}\|_F.
\end{equation}

Define $\gamma \doteq \polylog(n)$ that may vary from place to place and $\widehat{U}\doteq\big[ \widehat{\u}_{1}~ \widehat{\u}_{2}~\cdots~ \widehat{\u}_{r}\big]$. Denote $\underline{c}\doteq\min_{p\in[r]}{c_p^\star}^{1/3}$, $\bar{c}\doteq\max_{p\in[r]}{c_p^\star}^{1/3}$, and $\omega = \bar{c}/\underline{c}$. We are ready to introduce the assumptions needed to prove our main theorem as follows.

\begin{assumption}\label{asp:incoherence}
(Incoherence condition).	 The vector factors $\widehat{\u}$ in the target tensor $T^\star$ satisfy
\begin{align*}
\max_{i \neq j} | \langle\widehat{\u}_{i} ,\widehat{\u}_{j} \rangle |
\leq\frac{\gamma}{\sqrt{n}}.	
\end{align*}
\end{assumption}

\begin{assumption}\label{asp:spectrum}
(Bounded spectrum). The spectral norm of $\widehat{U}$ is bounded above as 
\begin{align*}
\|\widehat{U}\| \leq 1+c_1 \sqrt{\frac{r}{n}}.	
\end{align*} 	
\end{assumption}

\begin{assumption}\label{asp:isometry}
(Isometry of Gram-matrix). The Gram matrix satisfies the following isometry property  
\begin{align*}
\|(\widehat{U}^\top \widehat{U})\odot({\widehat{U}}^\top\widehat{U})-\eye_r\|\leq\frac{\gamma\sqrt r}{n}.
\end{align*} 	
where $\odot$ is the Hadamard product.
\end{assumption}

\begin{assumption}\label{asp:good_current}
(Warm start). The distance between the current variable $U$ and the matrix factor $U^\star$ is bounded with
\begin{align*}
\dist(U,U^\star)\leq
0.07\frac{m}{M}\frac{\underline{c}}{\omega^3}.
\end{align*} 
\end{assumption}

We remark that Assumptions \ref{asp:incoherence}-\ref{asp:isometry} hold with high probability if the factors $\{\widehat{\u}_p\}_{p=1}^r$ are generated independently according to the uniform distribution on the unit sphere~\cite[Lemmas 25, 31]{anandkumar2015learning}. 

\subsection{Main Results}

We now present our main theorem in the following: 
\begin{theorem} \label{thm:regularity}
Suppose that a $(r,m,M)$-restricted strongly convex and smooth function $f:\R^{n\times n\times n}\rightarrow \R$ has
a unique global minimizer at $T^\star$, which admits a CP decomposition $T^\star=U^\star\circ U^\star\circ U^\star \in \R^{n\times n\times n}$ as given in \eqref{eqn:truemodel}. Then, under Assumptions ~\ref{asp:incoherence}-\ref{asp:good_current} and in addition assuming $r = O\big(n^{1.25}\big)$, the following local regularity condition holds for sufficiently large $n$:
\begin{equation}
	\begin{aligned}
	\langle \nabla_U f(U\circ U\circ U), U-U^\star P_U\rangle 
	\geq & \frac{1}{2}\eta\|\nabla_U f(U\circ U\circ U)\|_F^2\\
	&+0.13m\underline{c}^4 \dist(U,U^\star)^2,	
	\end{aligned}
\end{equation}
as long as 
\begin{align}\label{eqn:stepsize:vary}
\eta\leq\frac{1}{18\|[\nabla f(T)]_{(1)}\|\cdot\|U\|+9M\|U\|^4}.
\end{align}
Here $T = U\circ U\circ U$ and $[\nabla f(T)]_{(1)}$ denotes the matricization of $\nabla f(T)$ along the first dimension.
\end{theorem}

The local regularity condition further implies linear convergence of the gradient descent algorithm~\eqref{eqn:algorithm} in a neighborhood of the true tensor factors $U^\star$ with proper choice of the stepsize, as summarized in the following two corollaries.

\begin{cor}[Linear Convergence with adaptive stepsize] 
Under the same assumptions as in Theorem~\ref{thm:regularity},  we have the following (adaptive) linear convergence  
\begin{equation}
\begin{aligned}\label{eqn:converge}
\dist(U^+,U^\star)^2
&\leq(1-0.26\eta m\underline{c}^4)\dist(U,U^\star)^2
\\
&\doteq \alpha(\eta)\cdot \dist(U,U^\star)^2
\end{aligned}
\end{equation}
when we run the gradient descent algorithm~\eqref{eqn:algorithm} with the stepsize $\eta$ satisfying \eqref{eqn:stepsize:vary}.  
\label{cor:1}
\end{cor}

\begin{cor}[Linear Convergence with constant stepsize]
Under the same assumptions as in Theorem~\ref{thm:regularity}, except that Assumption~\ref{asp:good_current} is replaced by a good initial condition:
\begin{equation}
\dist(U^0,U^\star)\leq 0.07\frac{m}{M}\frac{\underline{c}}{\omega^3},
\end{equation}
and the stepsize selection requirement~\eqref{eqn:stepsize:vary} is replaced with the constant stepsize satisfying
$\eta_0= \frac{1}{21.6M\|U^0\|^4},$
 the sequence $\{U^t: t=0,1,2,\cdots\}$ generated by
$$
U^{t+1}=U^{t}-\eta_0\nabla f(U^{t}\circ U^{t}\circ U^{t}),~~t=0,1,2,\cdots
$$
satisfies
\begin{align}\label{eqn:converge:2}
\dist(U^+,U^\star)^2
&\leq\alpha(\eta_0)\cdot\dist(U,U^\star)^2
\end{align}
with $\alpha(\eta_0)\doteq 1-0.26\eta_0 m\underline{c}^4$.
\label{cor:2}
\end{cor}

As a consequence, we conclude that solving the factored problem~\eqref{eqn:paramter} using the gradient descent algorithm~\eqref{eqn:algorithm} with a good initialization is guaranteed to recover the tensor factor matrix $U^\star$ with high probability if $r = O\big(n^{1.25}\big)$. The proof of the above theorem and corollaries can be found in supplementary material.

\section{Global Convergence}
\label{sec:global}

The local convergence analysis of applying vanilla gradient descent to  tensor optimization, though  developed for a class of sufficiently general problems, is not completely satisfactory as a good initialization might be difficult to find. Therefore, we are also interested in characterizing the global optimization landscape for these problems. Considering the difficulty of this task, we focus on a special case where the ground-truth third-order tensor admits an orthogonal decomposition and we are interested in finding its best rank-one approximation. We aim to characterize all its critical points and classify them into local minima, strict saddle points, and degenerate saddle points if there is any. We also want to exploit the properties of critical points to design a provable and efficient tensor decomposition algorithm.

\subsection{Main Results}
Consider the best rank-one approximation problem of an orthogonally decomposable tensor: 
\begin{align}
g(\u)=\|\u\otimes\u\otimes\u-T^\star\|_F^2,
\label{eqn:T}
\end{align}
 where $T^\star=\sum_{i=1}^r\u_i^\star\otimes\u_i^\star\otimes\u_i^\star$ and these true tensor factors $\{\u_i^\star\}$ are orthogonal to each other. This is a special case of \eqref{eqn:paramter} (and \eqref{eqn:tendec}). We characterize all possible critical points and their geometric properties in the following theorem: 
\begin{theorem}
Assume $T^\star=\sum_{i=1}^r\u_i^\star\otimes\u_i^\star\otimes\u_i^\star$, where $\{\u_i^\star\}$ are orthogonal to each other. Then any critical point $\widehat\u$ of $g(\u)$ in \eqref{eqn:T} takes the form $\widehat\u=\sum_{i=1}^r \lambda_i \u_i^\star$ for  $\vlambda\doteq[\lambda_1~\cdots~\lambda_r]^\top\in\R^r$ and
\begin{enumerate}
\item when $\|\vlambda\|_0=0$, $\widehat\u=\zero$ is a third-order saddle point, i.e., $\nabla^2 g(\widehat\u)=\zero$ and $\nabla^3 g(\widehat\u)\neq \zero$;

\item when $\|\vlambda\|_0=1$,  $\widehat\u=\u_{i}^\star$ with $i\in\{1,2,\ldots,r\}$ is a strict local minimum;

\item when $\|\vlambda\|_0\ge 2$, $\widehat\u$ is a strict saddle point, i.e., $\nabla^2 g(\u)$ has a negative eigenvalue.

\end{enumerate}
Here the $\ell_0$ ``norm" $\|\cdot\|_0$ counts the number of non-zero entries in a vector. Analytic expression for $\vlambda$ is given in the proof. 
\label{thm:2}
\end{theorem}

Theorem \ref{thm:2} implies that all second-order critical points are the true tensor factors except for zero. Based on this, we  develop a provable conceptual tensor decomposition algorithm as follows:
\begin{algorithm}[tb]
\caption{Iterative Gradient Descent for  Tensor Decomposition}
\label{alg:tensor}
\textbf{Input}: $T^\star$\\
\textbf{Initialization}: $T=T^\star$, $\widehat\u=\zero$\\
\textbf{Output}: Estimated factors $\{\u_i^\star\}$
\begin{algorithmic}[1] 
\STATE Let $i=0$.
\WHILE{${T}\neq\zero$}
\IF {$\widehat\u\neq\zero$}
\STATE $i=i+1.$
\STATE $\u_i^\star=\widehat\u$
\STATE ${T}\leftarrow{T}-\frac{\lg {T},\widehat\u\otimes\widehat\u\otimes\widehat\u\rg \widehat\u\otimes\widehat\u\otimes\widehat\u}{\|\widehat\u\|_2^3}$
\ENDIF
\STATE Find a second-order stationary point $\widehat\u$ of $g(\u)=\|\u\otimes\u\otimes\u-{T}\|_F^2$.
\ENDWHILE
\STATE \textbf{return} solution
\end{algorithmic}
\end{algorithm}

\begin{cor}
Assume $T^\star$ is a third-order orthogonal tensor with the tensor factors $\{\u_i^\star\}$.	Then with the input $T^\star$, \Cref{alg:tensor} almost surely recovers all the tensor factors $\{\u_i^\star\}$.
\label{cor:3}
\end{cor}
\begin{proof}[Proof of \Cref{cor:3}]
 It mainly follows from the many iterative algorithms can find a second-order stationary point \cite{lee2016gradient, li2019alternating, li2019provable, nesterov2006cubic, jin2017escape}. Then by \Cref{thm:2}, applying these iterative algorithms to $g(\u)$,  it converges to either a true tensor factor $\u_i^\star$ for $i\in[r]$ or the zero point (as a third-order saddle point is essentially a second-order stationary point). If it converges to a nonzero point, it must be a true tensor factor and we record it.  Then we can remove this component by projecting the target tensor $T$ into the orthogonal complement of $\u_i^\star$.  We repeat this process to the new deflated tensor until we get a zero deflated tensor. That means,  we have found all the true factors $\{\u_i^\star\}$.
\end{proof}

\subsection{Proof of \Cref{thm:2}}
 
Recall that 
\begin{equation*}
\begin{aligned}
T^\star=\sum_{i=1}^r\u_i^\star\otimes\u_i^\star\otimes\u_i^\star\doteq\sum_{i=1}^r\lambda_i\widehat\u_i\otimes\widehat\u_i\otimes\widehat\u_i.	
\end{aligned}	
\end{equation*}
Without loss of generality, we can extend the orthonormal set $\{\widehat\u_i\}_{i=1}^r$  to $\{\widehat\u_i\}_{i=1}^n$ as a full orthonormal basis of $\R^n$ and define 
\[\lambda_{i}\doteq0,~i\in[r]^c\doteq\{r+1,\ldots,n\}.\] Then, we have
$
T^\star=\sum_{i=1}^n\lambda_i\widehat\u_i\otimes\widehat\u_i\otimes\widehat\u_i.
$ 
Since $\{\widehat\u_i\}_{i=1}^n$ is a full orthonormal basis of $\R^n$, $\widehat U\doteq[\widehat\u_1\cdots~\widehat\u_n]$ is an orthonormal matrix, i.e., $\widehat U \widehat U^\top=\eye$. Then 
the best rank-1 tensor approximation problem is equivalent to
\begin{equation}
\begin{aligned}
g(\u)&=\|\u\otimes\u\otimes\u-T^\star\|_F^2
\\
&=\Big\|(\widehat U\widehat U^\top\u)\otimes(\widehat U \widehat U^\top\u)\otimes (\widehat U\widehat U^\top\u)
\\
&~~~~~~~~~-\sum_{i=1}^n\lambda_i\widehat \u_i\otimes\widehat \u_i\otimes\widehat \u_i\Big\|_F^2.	
\end{aligned}
\end{equation}

Expanding the squared norm and using the fact that $\widehat U$ is orthonormal, we get
\begin{equation}
\begin{aligned}
g(\u)&=\|(\widehat U^\top\u)\otimes(\widehat U^\top\u)\otimes(\widehat U^\top\u)-\diag_3(\vlambda)\|_F^2
\\
&\doteq\widehat g(\widehat U^\top \u)
\end{aligned}
\end{equation}
where we denote 
\begin{equation*}
\begin{aligned}
&\widehat g(\u)\doteq\|\u\otimes\u\otimes\u-\diag_3(\vlambda)\|_F^2	,
\\
&\diag_3(\vlambda)\doteq\sum_{i=1}^n \lambda_i\e_i\otimes\e_i\otimes\e_i	.
\end{aligned}
\end{equation*}

\begin{lem}\label{prop:eq}
The landscape of $g(\u)$ and $\widehat g(\u)$ are rotationally equivalent: $\u$ is a first/second-order stationary point of $g$ if and only if $\widehat U^\top \u$ is a first/second-order stationary point of $\widehat g$.

\end{lem}
\begin{proof}[Proof of \Cref{prop:eq}]
Since $g(\u) = \widehat g(\widehat U^\top\u)$,  by chain rule, 
\begin{equation}
\begin{aligned}
& \nabla  g(\u)= \widehat U \nabla \widehat g(\widehat U^\top \u),
\\
&\nabla^2 g(\u)= \widehat U \nabla^2 \widehat g(\widehat U^\top \u) \widehat U^\top.
\end{aligned}
\end{equation}
 Then it directly follows from the definitions of first/second stationary points. \end{proof}

Therefore by \Cref{prop:eq},  to understand the landscape of $g(\u)$, it suffices to study that of 
\begin{equation*}
\widehat g(\u)=\|\u\otimes\u\otimes\u-\diag_3(\vlambda)\|_F^2.
\end{equation*}
We compute its derivatives up to third-order:
\begin{equation*}
\begin{aligned}
\nabla \widehat g(\u)=&6\|\u\|_2^4\u- 6\vlambda\odot\u\odot\u,\\
\nabla^2\widehat g(\u)=&
6\|\u\|_2^4\eye+24\|\u\|_2^2\u\u^\top-12\diag_3(\vlambda\odot\u),
\\
\nabla^3 \widehat g(\u)=& 24\|\u\|_2^2 \tmop{Sym}(\eye\otimes\u)+48\u\otimes\u\otimes\u-12\diag_3(\vlambda),	
\end{aligned}
\end{equation*}
where $\odot$ is the Hadamard product and $\tmop{Sym}(T)$ is the sum of all the three permutations of $T$.

Now define $J$ as the index set of any critical point $\u$ such that $u_i\neq 0$ for $i\in J$, i.e.,
\[\u_J\neq\zero,~\u_{J^c}=\zero,~\|\u\|_0=|J|.\]
By the critical point equation 
\begin{equation}
\begin{aligned}
\widehat\u\|_2^4\widehat\u-\vlambda\odot \widehat\u\odot\widehat\u = 	\zero
\end{aligned}	
\end{equation}
and $\lambda_i=0$ for $i\in[r]^c$, we conclude that $J\subset[r]$. In the following, we divide the problem into three cases: $|J|=0,|J|=1,~\text{and}~|J|\ge2$.

\begin{itemize}
\item Case I: $|J|=0$. That is $\widehat\u=\zero$. Then, we have 
\begin{equation*}
\begin{aligned} 
&\nabla \widehat g(\widehat\u)=\zero,
~~\tmop{and}~~\nabla^2 \widehat g(\widehat\u)=\zero,
\end{aligned}	
\end{equation*}
but 
\begin{equation*}
\begin{aligned} 
\nabla^3\widehat g(\widehat\u)=-12\diag_3(\vlambda). 
\end{aligned}	
\end{equation*}
This implies $\widehat\u=\zero$ is a third-order saddle point of $\widehat g$.
\item Case II: $|J|=1$. Since $J\subset[r]$, 
 let $J=\{k\}$ for some $k\in[r]$. Then,
 \begin{equation*}
\begin{aligned} 
&\widehat\u=\widehat u_k \e_k~~\tmop{and}~~\nabla \widehat g(\widehat\u)=6\widehat u_k ^5\e_k -6\lambda_k \widehat u_k ^2\e_k =\zero,
\end{aligned}	
\end{equation*}
which implies that $\widehat u_k =\sqrt[3]{\lambda_k }$. We also have 
 \begin{equation*}
\begin{aligned}
\nabla^2 \widehat g(\widehat\u)&=6\lambda_k^{4/3}\eye+24 \lambda_k^{4/3} \e_k \e_k^\top-12\lambda_k^{4/3} \e_k \e_k ^\top\\
&=6\lambda_k^{4/3}\eye+12 \lambda_k^{4/3} \e_k \e_k^\top\succ0
\end{aligned}	
\end{equation*}
Therefore, any critical point $\widehat\u$ with $\|\widehat\u\|_0=1$ has the form  $\widehat\u=\sqrt[3]{\lambda_k } \e_k$ for $k\in[r]$, and is a strict local minimum of $\widehat g$. 

\item Case III: $|J|\ge 2$. Also, we know that $J\subset[r]$. With the critical point equation, we get 
 \begin{equation*}
\begin{aligned} 
\vlambda_J\odot\widehat\u_J=\|\widehat\u\|_2^4 \one_{|J|}.
\end{aligned}	
\end{equation*}
Further notice that
 \begin{equation*}
\begin{aligned} 
\diag_3(\vlambda_J\odot\widehat\u_J)=\|\widehat\u\|_2^4\eye_{|J|}.
\end{aligned}	
\end{equation*}
Plugging this to the sub-Hessian 
 \begin{equation*} 
[\nabla^2 \widehat g(\widehat\u)]_{J,J}=24\|\widehat\u\|_2^2\widehat\u_J\widehat\u_J^\top-6\|\widehat\u\|_2^4\eye_{|J|}.
\end{equation*}
Now for any $\d\in\R^n$ with $\d_J^\top\u_J=0, \d_{J^c}=\zero$, we have 
 \begin{equation*}
[\nabla^2 \widehat{g}(\widehat\u)](\d,\d)=-6\|\widehat\u\|_2^4\|\d_J\|_2^2=-6\|\widehat\u\|_2^4\|\d\|_2^2,
\end{equation*}
 implying that 
 \begin{equation*}
 \lambda_{\min}(\nabla^2 \widehat g(\widehat\u))\le-6\|\widehat\u\|_2^4<0.
 \end{equation*}
  Therefore, any critical point $\widehat\u$ with $\|\widehat\u\|_0\ge 2$ has the form $\widehat\u_J=\frac{\|\widehat\u_J\|_2^4}{\vlambda_J}$ (pointwise), and is a strict saddle point of $\widehat g$.
\end{itemize}
Together with \Cref{prop:eq}, we complete the proof of \Cref{thm:2}.

\section{Numerical Experiments}
\label{sec:simu}

\paragraph{Computing Infrastructure} All the numerical experiments
are performed on a 2018 MacBook Pro with operating system of macOS version 10.15.7, processor of 2.6 GHz 6-Core Intel Core i7, memory of 32 GB, and MATLAB version of R2020a.

\begin{figure*}[!htp]
\centering	
\includegraphics[width=1\textwidth]{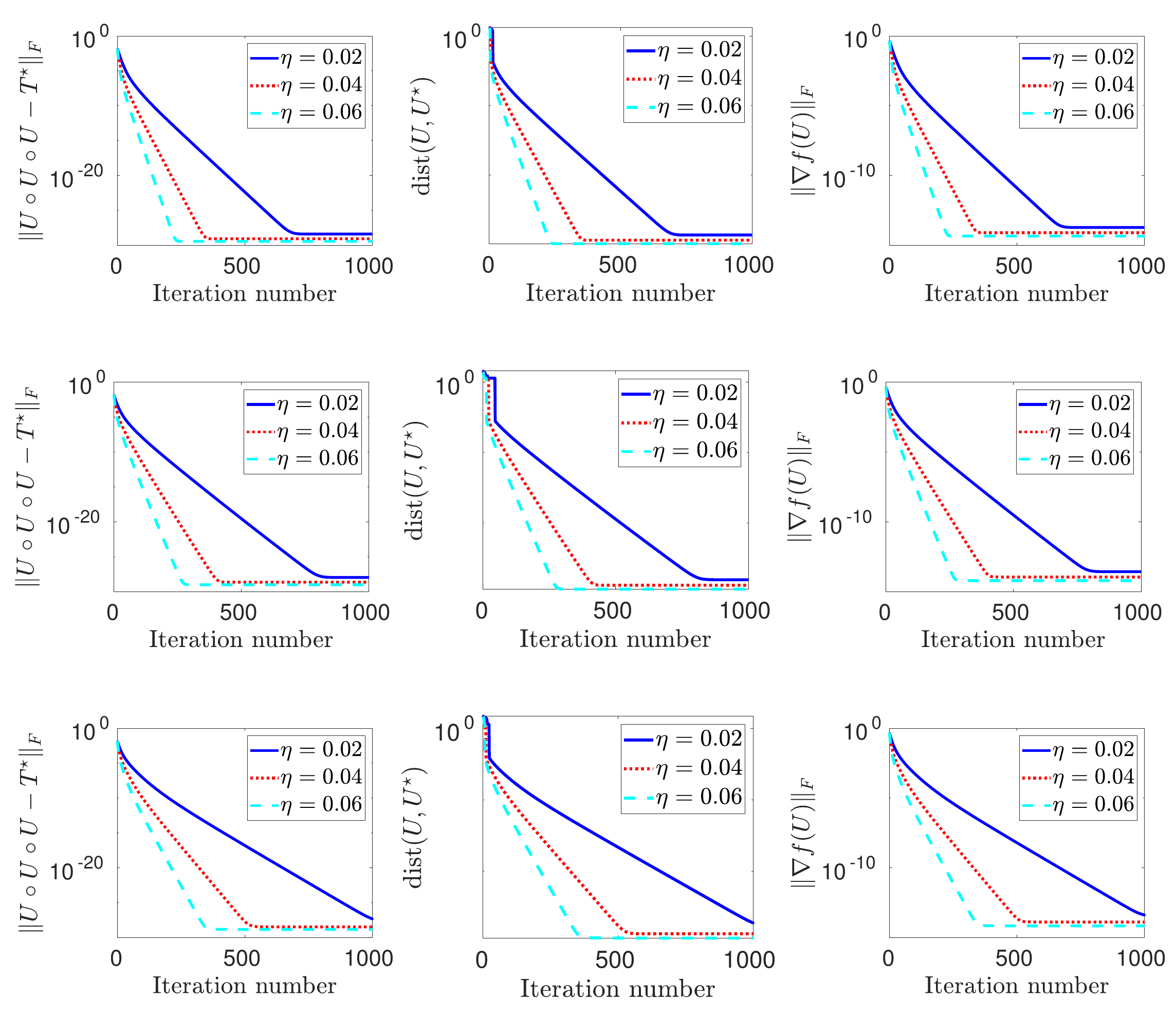}
\caption{Linear convergence of gradient descent when applied to tensor factorization problem~\eqref{eqn:tendec}. Here, $r=n/2$ ({\em top row}),   $r=n$ ({\em middle row}), and  $r=3n/2$ ({\em bottom row}) with $n=64$.  We initialize the starting point as $U^\star + \alpha D$ with $\alpha=0.07$ and set $D$ as a normalized Gaussian matrix with $\|D\|_F=1$. We record the three metrics $\|\nabla f(U)\|_F$ ({\em left column}), $\|U\circ U \circ U - T^\star\|_F$ ({\em middle column}), and $\dist(U,U^\star)$ ({\em right column}) for total $10^3$ iterations with different stepsize $\eta$,  which is consistent with the linear convergence analysis of gradient descent on general Burer-Monteiro tensor optimizations. 
}
\label{Fig}
\end{figure*}

In the first experiment, we illustrate the linear convergence of the gradient descent algorithm within the contraction region $\dist(U^0,U^\star)\le 0.07\frac{m}{M}\frac{\underline{c}}{\omega^3}$ in solving the tensor decomposition problem~\eqref{eqn:tendec}, where $M=m=1$ in this case. 
We set $n = 64$ and vary $r$ with three different values: $n/2,~n,~3n/2$ to get an undercomplete, complete, and overcomplete target tensor $T^\star$, respectively.  
We generate the $r$ columns of $U^\star$ independently according to the uniform distribution on the unit sphere and form $T^\star=U^\star\circ U^\star\circ U^\star$. According to \cite[Lemmas 25, 31]{anandkumar2015learning} and \Cref{cor:2}, if $\dist(U^0,U^\star)\le 0.07\frac{m}{M}\frac{\underline{c}}{\omega^3}=0.07$ (because $\|\u_i^\star\|_2=1$ implies $\bar{c}=\underline{c}=\omega=1$), the gradient descent with a sufficiently small constant stepsize would converge linearly to the true factor $U^\star$. To illustrate this,  we initialize the starting point as $U^\star + \alpha D$ with $\alpha=0.07$ and set $D$ as a normalized Gaussian matrix with $\|D\|_F=1$. We record the three metrics $\|\nabla f(U)\|_F$, $\|U\circ U \circ U - T^\star\|_F$, and $\dist(U,U^\star)$ for total $10^3$ iterations with different stepsizes $\eta$ in Figure~\ref{Fig}, which is consistent with the linear convergence analysis of gradient descent on general Burer-Monteiro tensor optimizations  in \Cref{cor:2}. 

In the second experiment, with the same settings as above except varying $\alpha$, 
we record the success rate by running 100 trials for each fixed $(r,\alpha)$-pair and declare one successful instance if the final iterate $U$ satisfies $\dist(U,U^\star)\le 10^{-3}$. We repeat these experiments for different  $\alpha\in\{0.5,1,2,4,8,16\}$.
\Cref{table:1} shows that when $\alpha$ is small enough ($\alpha\le 2$),  the success rate is 100\% for all the undercomplete ($r=n/2$), complete ($r=n$), and overcomplete ($r=3n/2$) cases; and when $\alpha$ is comparatively large ($\alpha\in[4,8]$),  the success rate degrades dramatically when $r$ increases. Finally, when $\alpha$ is larger than certain threshold, the success rate is 0\%. This in consistence with Corollaries \ref{cor:1} and \ref{cor:2}.

\begin{table}[!htp]
\def\ROWCOLOR{black!10!white}
\centering
\caption{Success ratio with $\eta=0.02$ ({\em top}),   $\eta=0.04$ ({\em middle}), and  $\eta=0.06$ ({\em bottom}).} 
\begin{tabular}{l c c c c c c} 
 \hline
   $\alpha$ &     0.5 & 1 & 2 & 4 & 8 & 16\\  
 \hline\hline
  $r=n/2$ &     {100\%} & 100\% & 100\%  & 100\%  &	   100\%  &	    {0\%}\\ 
  $r=n$ &     {100\%} & 100\% & 100\% & 100\%  &100\%  &	    {0\%}\\
   $r=3n/2$ &     {100\%} & 100\% & 100\% & 100\%  &	   5\%  &	    {0\%}\\ 
 \hline 
\end{tabular}
\\[2ex]
\centering
\begin{tabular}{l c c c c c c} 
 \hline
   $\alpha$ &     0.5 & 1 & 2 & 4 & 8 & 16\\
 \hline\hline
  $r=n/2$ &     {100\%} & 100\% & 100\%  & 100\%  &	   100\%  &	    {0\%}\\ 
   $r=n$ &     {100\%} & 100\% & 100\% & 100\%  &	   0\%  &	    {0\%}\\
  $r=3n/2$ &     {100\%} & 100\% & 100\% & 100\%  &	   0\%  &	    {0\%}\\     
 \hline
\end{tabular}
\\[2ex]
\begin{tabular}{l c c c c c c} 
 \hline
   $\alpha$ &     0.5 & 1 & 2 & 4 &8 & 16\\  
 \hline \hline
   $r=n/2$ &     {100\%} & 100\% & 100\%  & 100\%  &~ 38\%  & {0\%}\\ 
  $r=n$ &     {100\%} & 100\% & 100\% & 100\%  &~ 0\%  & {0\%}\\
  $r=3n/2$ &     {100\%} & 100\% & 100\% & 83\%  &~ 0\%  &{0\%}\\   
 \hline
\end{tabular}
\label{table:1}
\end{table}

\section{Conclusion}
\label{sec:conc}

In this work, we investigated the local convergence of third-order tensor optimization with general convex and well-conditioned objective functions. Under certain incoherent conditions, we proved the local regularity condition for the nonconvex factored tensor optimization resulted from the Burer-Monteiro reparameterization. We highlighted that these assumptions are satisfied for randomly generated tensor factors. With this local regularity condition, we further provided a linear convergence analysis for the gradient descent algorithm started in a neighborhood of the true tensor factors. Complimentary to the local analysis, we also presented a complete characterization of the global optimization landscape of the best rank-one tensor approximation problem.

\section*{Acknowledgments}
S. Li  gratefully acknowledges support from the National Science Foundation (NSF award DMS 2011140). We thank Prof. Gongguo Tang (CU Boulder) for fruitful discussions.

\bigskip

\bibliography{nips_tensor}

\newpage
\appendix

\setcounter{secnumdepth}{2}

\onecolumn

\section{Proof of \Cref{thm:regularity}}

For convenience, we present all the necessary tensor/matrix notations in \Cref{table:notation}, where we assume $U=[\u_1~\u_2~\cdots~\u_r]\in\R^{n\times r}$, $V=[\v_1~\v_2~\cdots~\v_r]\in\R^{m\times r}$, and $W=[\w_1~\w_2~\cdots~\w_r]\in\R^{l\times r}$.
\begin{table}[ht]
\caption{Tensor/matrix notations}
\centering
\begin{tabular}{c l l}
\hline\hline
Symbols & Meaning  & Explanation \\ [0.5ex] 
\hline \\[-2ex]
$\otimes$&outer/tensor product&   $\u\otimes\v\otimes\w\doteq T\in\R^{n\times m\times l}$ with $T_{ijk}=u_iv_jw_k$.  \\ \\
$\circ$&group tensor product&$U\circ V=UV^\top$; $U\circ V \circ W\doteq\sum_{j=1}^r \u_j\otimes \v_j\otimes\w_j$ .\\ \\
$\bigotimes$& Kronecker products& $U\bigotimes V\doteq\begin{bmatrix} u_{11}V& \cdots & u_{1r}V\\
\vdots& \ddots & \vdots\\
u_{n}V& \cdots & u_{n r}V
\end{bmatrix}\in\R^{n m\times r^2}$. \\ \\
$\circledast$&Khatri-Rao product& $U \circledast V\doteq\begin{bmatrix}\u_1\bigotimes \v_1&\cdots&\u_m\bigotimes\v_m\end{bmatrix}\in\R^{nm\times r}$ . \\ \\
$\odot$& Hadamard product& $U \odot V\doteq M$ with $M_{ij}=u_{ij}v_{ij}$ when $m=n$.
\\
[1ex]
\hline
\end{tabular}
\label{table:notation}
\end{table}

 To prove \Cref{thm:regularity}, we need the following key lemmas.\footnote{Lemmas~\ref{lem:TTs} and \ref{lem:TTs_inner} are proved will be proved in the remaining sections. Lemma \ref{lem:FgradT} follows from the fact that $\nabla f(T^\star)=\mathbf{0}$, $\rank(T^\star)\le r$, and the $r$-restricted  $M$-smoothness. Lemma \ref{lem:ABC_inner} is a result from the {\em H\"{o}der's in equality} and the fact that for any 
 tensor $T=\sum_{j\in[r]}\a_j\otimes\b_j\otimes \b_j$, its tensor nuclear norm can be bounded as $\|T\|_*\le \sum_{j\in[r]}\|\a_j\|_2\|\b_j\|_2^2\le \left(\max_{j\in[r]}\|\a_j\|_2\right)\cdot \|B\|_F^2$ by the definition of tensor nuclear norm \cite{li2016super}.}
\begin{lem}\label{lem:TTs}
Denote $H\doteq U-U^\star P_U$ and assume $\|H\|_F\leq0.07\frac{\underline{c}}{\omega^3}$.\footnote{Note that we use a weaker assumption $\|H\|_F\leq0.07\frac{\underline{c}}{\omega^3}$ here when compared with the one used in \Cref{thm:regularity}, i.e., $\|H\|_F\leq0.07 \frac{m}{M}\frac{\underline{c}}{\omega^3}$.} Under Assumptions \ref{asp:incoherence}-\ref{asp:good_current}, if $r=O(n^{1.25})$, we have the following bounds hold for sufficiently large $n$:
$$
1.679\underline{c}^4\|H\|_F^2 \leq \|T-T^\star\|_F^2\leq  10.336\omega^4\underline{c}^4\|H\|_F^2.
$$
\end{lem}

\begin{lem}\label{lem:TTs_inner}
Suppose that a $(r,m,M)$-restricted strongly convex and smooth function $f:\R^{n\times n\times n}\rightarrow \R$ has
a unique global minimizer at $T^\star$ of rank at most $r$. Then for any  $\eta\leq\frac{1}{18\|[\nabla f(T)]_{(1)} \|\|U\|+9M\|U\|^4},$  we have\footnote{Note that $\nabla f(T) = \nabla f(U\circ U\circ U)$ is the gradient of $f$ with respect to a tensor $T = U\circ U\circ U$, while $\nabla_U f(U\circ U\circ U)$ is the gradient of $f$ with respect to a matrix $U$.}
$$
\langle \nabla f(U\circ U\circ U) ,U\circ U\circ U-T^\star\rangle \geq \frac{1}{2}\eta\|\nabla_U f(U\circ U\circ U)\|_F^2+\frac{m}{2}\|U\circ U\circ U-T^\star\|^2_F.
$$
\end{lem}

\begin{lem}\label{lem:FgradT}
Suppose that a $(r,m,M)$-restricted strongly convex and smooth function $f:\R^{n\times n\times n}\rightarrow \R$ has
a unique global minimizer at $T^\star$ of rank at most $r$. Then for any $T$ with $\rank(T)\le r$,  we have  
$$
\|\nabla f(T)\|_F=
\|\nabla f(T)-\nabla f(T^\star)\|_F
\leq M\|T-T^\star\|_F.
$$
\end{lem}

\begin{lem}\label{lem:ABC_inner}
For any two matrices $A=[\a_1~\cdots~\a_r]\in\R^{n\times r}$, $B=[\b_1~\cdots~\b_r]\in\R^{n\times r}$ and a tensor $C\in\R^{n\times n \times n}$, we have the following inequality holds
\begin{align*}
|\langle C, A\circ B\circ B\rangle|
\leq  \|C\|\cdot \max_{j\in[r]} \|\a_j\|_2 \cdot \|B\|_F^2.
\end{align*}
\end{lem}

\begin{proof}[Proof of \Cref{thm:regularity}]

Observe that 
\begin{align*}
\langle \nabla_U f, U-U^\star P_U\rangle
=&\langle 3[\nabla f]_{(1)}(T) U\circledast U, U-U^\star P_U\rangle
\\
=&3\langle[\nabla f]_{(1)}(T) , U(U\circledast U)^T-U^\star P_U(U\circledast U)^T\rangle
\\
=&3\langle \nabla f(T) ,T-U^\star P_U \circ U\circ U\rangle
\text{ (Denote new $U^\star=U^\star P_U.$)}
\\
=&\langle \nabla f(T),(T-T^\star)+2(U-U^\star)\circ (U-U^\star)\circ (U-U^\star)+3U^\star\circ (U-U^\star)\circ(U-U^\star)\rangle
\\
=&\langle \nabla f(T),(T-T^\star)+2H\circ H\circ H+3U^\star\circ H\circ H\rangle
\\
=&\langle \nabla f(T),T-T^\star\rangle
+\langle \nabla f(T),2H\circ H\circ H\rangle
+\langle \nabla f(T),3U^\star\circ H\circ H\rangle.
\end{align*}
Next, we bound the above three terms in sequence. Using Lemma \ref{lem:TTs_inner}, the first term can be bounded with
$$
\langle \nabla f(T) ,T-T^\star\rangle \geq \frac{1}{2}\eta\|\nabla_U f\|_F^2+\frac{m}{2}\|T-T^\star\|^2_F.
$$
Combining Lemmas \ref{lem:FgradT} and \ref{lem:ABC_inner}, the absolute value of the last two terms can be bounded with
\begin{align*}
|\langle \nabla f(T),2H\circ H\circ H\rangle|
\leq &2\|\nabla f(T)\|\|H\|_F^3 
\\
\leq & 2M\|T-T^\star\|_F\|H\|_F^3,
\\ 
\\
|\langle \nabla f(T),3U^\star\circ H\circ H\rangle|
\leq  &3\bar{c}\|\nabla f(T)\|\|H\|_F^2 
\\
\leq &3\bar{c}M\|T-T^\star\|_F\|H\|_F^2.
\end{align*}
Then, we have 
\begin{align*}
\langle \nabla_U f, U - U^\star P_U\rangle   
&\geq \frac{1}{2}\eta\|\nabla_U f\|_F^2 + \bigg( \frac{m}{2}\|T - T^\star \|_F^2 - 2M\|T - T^\star \|_F\|H\|_F^3 - 3M\bar{c}\|T - T^\star \|_F\|H\|_F^2 \bigg)\\
&\doteq \frac{1}{2}\eta\|\nabla_U f\|_F^2+\Pi.    
\end{align*}

Applying Lemma~\ref{lem:TTs}, we obtain
\begin{align*}
\Pi&\geq
\frac{1}{2}
\bigg(
1.679m\underline{c}^2
-4\sqrt{10.336}M\omega^2\|H\|_F^2
-6\sqrt{10.336}M\omega^3\underline{c}\|H\|_F
\bigg)\underline{c}^2\|H\|_F^2\\
& \geq \bigg(1.679-2\sqrt{10.336}\bigg(\frac{2(0.07)^2}{\omega^4}\cdot\frac{m}{M} +3(0.07)\bigg)\bigg)\frac{1}{2}m\underline{c}^4\|H\|_F^2\\
& \geq \bigg(1.679-2\sqrt{10.336}(2(0.07)^2+3(0.07))\bigg)\frac{1}{2}m\underline{c}^4\|H\|_F^2\\
& \geq 0.13 m\underline{c}^4\|H\|_F^2,
\end{align*}
where the second inequality follows from $\|H\|_F\leq0.07\frac{\underline{c}}{\omega^3}$.
Finally, we can get the local regularity condition
$$
\langle \nabla_U f, U-U^\star P_U\rangle \geq \frac{1}{2}\eta\|\nabla_U f\|_F^2+0.13m\underline{c}^4 \dist(U,U^\star)^2
$$
and complete the proof of \Cref{thm:regularity}.
\end{proof}

\section{Proof of Corollaries \ref{cor:1} and \ref{cor:2}}
\label{sec:proof_cor_1}

We first prove \Cref{cor:1} by using the results stated in \Cref{thm:regularity}. In particular, we have
\begin{align*}
\dist(U^+,U^\star)^2
\leq& \|U^+ - U^\star P_U\|_F^2
\\
= &\|U^+-U+U-U^\star P_U\|_F^2\\
=& \|U^+-U\|_F^2+ \|U-U^\star P_U\|_F^2+2\langle U^+-U,U-U^\star P_U\rangle\\
=& \eta ^2\|\nabla_U f(U\circ U\circ U)\|_F^2+ \dist(U,U^\star)^2-2\eta \langle \nabla_U f(U\circ U\circ U),U-U^\star P_U\rangle\\
\leq & \eta ^2\|\nabla_U f(U\circ U\circ U)\|_F^2  +  \dist(U,U^\star)^2 - 2\eta \bigg( \frac{1}{2}\eta\|\nabla_U f(U\circ U\circ U)\|_F^2 + 0.13m\underline{c}^4 \dist(U,U^\star)^2 \bigg)\\
\leq &(1-0.26\eta m\underline{c}^4)\dist(U,U^\star)^2.
\end{align*}

Next, we continue to prove \Cref{cor:2}. With \Cref{cor:1}, it suffices to show that there is a stepsize $\eta_0$ such that
\begin{align}
\eta_0\leq \min\left\{\frac{1}{18\|[\nabla f(U^t\circ U^t\circ U^t)]_{(1)} \|\cdot\|U^t\|+9M\|U^t\|^4}\doteq s(U^t):~\forall t\ge 0 \right\}.
\label{eqn:QW}
\end{align}

For this purpose, we need to find an upper bound for $\|U^t\|$ and $\|[\nabla f(U^t\circ U^t\circ U^t)]_{(1)} \|$.

We first bound $\|U^t\|$. Note that $U^0$ lies in the local region if $\eta$ is small enough. Then, for any current variable $U^t$, we have
\begin{align}
\max\{\|U^0-U^\star\|_F,\|U^t-U^\star\|_F\}\leq0.07\frac{m}{M}\frac{\underline{c}}{\omega^3}\leq0.07\underline{c}\leq 0.07\|U^\star\|,    
\label{eqn:U0U}
\end{align}
which implies that
\begin{align*}
\max\{\|U^0-U^\star\|,\|U^t-U^\star\|\}\leq0.07\|U^\star\|,    
\end{align*}
and 
\begin{align*}
0.93\|U^\star\|\leq \max\{\|U^t\|,\|U^0\|\}\leq 1.07 \|U^\star\|.    
\end{align*}
Thus, we can get
\begin{align}
\|U^t\|\leq \frac{1.07}{0.93}\|U^0\| \text{ and } \bar{c}\leq\|U^\star\|\leq \frac{1}{0.93}\|U^0\|,
\label{eqn:UUs}    
\end{align}
where the first part follows from $\|U\|\le1.07\|U^\star\|$ and $\|U^\star\|\le \frac{1}{0.93} \|U^0\|$, and the second line follows from  $\|U^\star\|\le \frac{1}{0.93} \|U^0\|$, $\|U^\star\|\ge \max_i \|\u_i^\star\|_2$, and $\bar{c}\doteq\max_i \|\u_i^\star\|_2$.

Next, we bound $\|[\nabla f(U^t\circ U^t\circ U^t)]_{(1)} \|$, or equivalently $\|[\nabla f(T)]_{(1)} \|$ with $T=U^t\circ U^t\circ U^t$. Observe that
\begin{align}
\|[\nabla f(T)]_{(1)} \|
&\leq\|\nabla f(T)\|_F
\\
&\leq M\|T-T^\star\| \nonumber\\
&\leq(10.336)^{1/2} M \bar{c}^2\|U^t-U^\star\|_F \nonumber\\
&\leq (10.336)^{1/2} M\bar{c}^2\bigg(\frac{0.07}{0.93}\|U^0\|\bigg)\nonumber\\
&\leq 0.242 M\bar{c}^2\|U^0\|,\label{eqn:Nablas}
\end{align}
where the second inequality follows from Lemma~\ref{lem:FgradT}, the third inequality follows from Lemma~\ref{lem:TTs}, and the fourth inequality follows from \eqref{eqn:U0U} and \eqref{eqn:UUs}.

Finally,  plugging the upper bound of $\|U^t\|$  in \eqref{eqn:UUs}  and the upper bound of  $\|[\nabla f(T)]_{(1)} \|$ in  \eqref{eqn:Nablas} to the definition of $s(U^t)$, i.e.,
\[s(U^t)\doteq\frac{1}{18\|[\nabla f(U^t\circ U^t\circ U^t)]_{(1)} \|\cdot\|U^t\|+9M\|U^t\|^4},\]
we can bound $s(U^t)$ as
\begin{align*}
s(U^t)
&\geq \frac{1}
{18(0.242 M(\frac{1}{0.93}\|U^0\|)^2\|U^0\|)(\frac{1.07}{0.93}\|U^0\|)+9M(\frac{1.07}{0.93}\|U^0\|)^4}\geq \frac{1}{21.6M\|U^0\|^4}.
\end{align*}
Therefore, we complete the proof of \Cref{cor:2}.

\section{Proof of Lemmas~\ref{lem:TTs} and \ref{lem:TTs_inner}}
\label{sec:proof_key_lems}

To prove Lemmas~\ref{lem:TTs} and \ref{lem:TTs_inner}, we require some additional auxiliary lemmas.

\subsection{Auxiliary lemmas}
\label{sec:auxi}

We first introduce the auxiliary lemmas used for proving \Cref{lem:TTs}, which are mainly related with 
Assumptions \ref{asp:incoherence}-\ref{asp:isometry}. Since Assumptions \ref{asp:incoherence}-\ref{asp:isometry}  are concerned with the normalized factors $\widehat{U}$, we need transform these assumptions to apply to $U^\star$ by using that the dynamic range of the coefficients $\{c_j^{\star1/3}\}$ is small.

\begin{lem}\label{lem:aux:incoherence}
Under  \Cref{asp:incoherence} on $\widehat U$, the mutual incoherence coefficient of $U^\star$ is upper bounded  by
\begin{align}
\mu^\star=\max_{i \neq j} | \langle{\u}_{i}^{\star} ,{\u}_{j}^{\star} \rangle |
\leq\frac{\bar{c}^{2} \gamma}{\sqrt{n}}.
\end{align}
\end{lem}

\begin{proof}[Proof of \Cref{lem:aux:incoherence}]
By \Cref{asp:incoherence}, we have
$$
\max_{i \neq j} | \langle\widehat{\u}_{i} ,\widehat{\u}_{j} \rangle |\leq\frac{\gamma}{\sqrt{n}}.
$$
Since $\u^\star_i={c_i^\star}^{1/3}\widehat{\u}_i$ for any $i\neq j$, we have
\begin{align*}
|\langle\u^\star_i ,\u^\star_j \rangle|
=&{c_i^\star}^{1/3}{c_j^\star}^{1/3}|\langle\widehat{\u}_{i} ,\widehat{\u}_{j} \rangle|
\\
\leq& \big(\max_{k\in[r]} {c_k^\star}^{1/3}\big)^2\cdot \max_{i \neq j} | \langle\widehat{\u}_{i} ,\widehat{\u}_{j} \rangle |
\\
=& \bar{c}^2\cdot\frac{\gamma}{\sqrt{n}}
\end{align*}
\end{proof}

\begin{lem}\label{lem:aux:spectrum}
Under  \Cref{asp:spectrum} on $\widehat U$, the operator norm of $U^\star$  is bounded by
\begin{align}
\|U^\star\|\leq \bar{c}\bigg(1+c_1 \sqrt{\frac{r}{n}}\bigg).
\end{align}
\end{lem}

\begin{proof}[Proof of \Cref{lem:aux:spectrum}]
By \Cref{asp:spectrum} on $\widehat{U}$, we have
$$
\|\widehat{U}\|\leq 1+c_1 \sqrt{\frac{r}{n}}.
$$
This along  with $\u^\star_i=c_i^{1/3}\widehat{\u}_i$ implies that
\begin{align*}
U^\star=\widehat{U}
\underbrace{\begin{bmatrix}
{c_1^\star}^{1/3}&0&\cdots&0\\
0&{c_2^\star}^{1/3}&\cdots&0\\
\vdots&\vdots&\vdots&\vdots\\
0&\cdots&0&{c_r^\star}^{1/3}
\end{bmatrix}}_{C}.
\end{align*}
Therefore,
\begin{align*}
\|U^\star\|
&\leq\|\widehat{U}\|\cdot\|C\|
\\
&\leq(1+c_1 \sqrt{\frac{r}{n}})\cdot\max_{k\in[r]} {c_k^\star}^{1/3}
\\
&=(1+c_1 \sqrt{\frac{r}{n}})\cdot\bar{c}.
\end{align*}
The second inequality follows from  \Cref{asp:spectrum} on $\widehat{U}$, and the fact that the operator norm of a diagonal matrix is the largest absolute value of its diagonal entries.
\end{proof}

Further, we claim that the $2\rightarrow 3$ norm and $2\rightarrow 4$ norm of $U^\star$
are also upper bounded under \Cref{asp:incoherence}, where the general $p\to q$ norm is denoted and defined by  $\|A\|_{p\to q}\doteq \max_{\x\neq \zero} \frac{\|A\x\|_q}{\|\x\|_p}. $

In particular, we have the following lemma.

\begin{lem}\label{lem:aux:p:norm}
Under  \Cref{asp:incoherence} on $\widehat U$,
the $2\rightarrow 3$ norm and $2\rightarrow 4$ norm of $U^\star$ are bounded with
\begin{align*}
&\|{U^\star}^\top\|^k_{2\rightarrow 3}
\leq \bar{c}^k\big(1+O(\gamma n^{-\epsilon})\big), \quad \tmop{if}~r=O(n^{1.25-1.5\epsilon}),
\\
&\|{U^\star}^\top\|^k_{2\rightarrow 4}
\leq \bar{c}^k\big(1+O(\gamma n^{-\epsilon})\big), \quad \tmop{if}~r=O( n^{1.5-2\epsilon}).
\end{align*}
\end{lem}

\begin{proof}[Proof of \Cref{lem:aux:p:norm}]
It suffices to show
\begin{align*}
\|{\widehat{U}}^\top\|_{2\rightarrow 3}
&\leq 1+O(\gamma n^{-\epsilon}), \quad \tmop{if}~r=O(n^{1.25-1.5\epsilon}),
\\
\|{\widehat{U}}^\top\|_{2\rightarrow 4}
&\leq 1+O(\gamma n^{-\epsilon}), \quad \tmop{if}~r=O( n^{1.5-2\epsilon}).
\end{align*}

First of all, identify that
\begin{align*}
\|{U^\star}^\top\|_{2\rightarrow 3}
&=\|\widehat{U}C\|_{2\rightarrow 3}
\\
&\leq \|\widehat{U}\|_{2\rightarrow 3}\cdot \|C\|
\\
&\leq \big(1+O(\gamma n^{-\epsilon})\big)\cdot  \bar{c},
\end{align*}
where the second inequality is because
\begin{align*}
\|AB\|_{2\rightarrow 3}
&=\max_{\|\x\|_2=1}\|AB\x\|_3
\\
&\leq  \|A\|_{2\rightarrow 3}\cdot\max_{\|\x\|_2=1}\|Bx\|_2
\\
&= \|A\|_{2\rightarrow 3}\cdot\|B\|.
\end{align*}
Since
\[
\big(1+O(\gamma n^{-\epsilon})\big)^k=1+O(\gamma n^{-\epsilon}),
\]
we have
$$
\|{U^\star}^\top\|^k_{2\rightarrow 3}\leq \bar{c}^k\big(1+O(\gamma n^{-\epsilon})\big).
$$
A similar argument applies to $\|{U^\star}^\top\|^k_{2\rightarrow 4}\leq \bar{c}^k\big(1+O(\gamma n^{-\epsilon})\big)$.

Therefore, to complete this proof, it suffices to show
\begin{align*}
\|{\widehat{U}}^\top\|_{2\rightarrow 3}
&\leq 1+O(\gamma n^{-\epsilon}), \quad \tmop{if}~r=O(n^{1.25-1.5\epsilon}), 
\\
\|{\widehat{U}}^\top\|_{2\rightarrow 4}
&\leq 1+O(\gamma n^{-\epsilon}), \quad \tmop{if}~r=O( n^{1.5-2\epsilon}),
\end{align*}
which are clearly true when $\|{\widehat{U}}^\top\|_{2\rightarrow p} \leq 1$ with $p=3,4$.

We continue to consider the case when $\|{\widehat{U}}^\top\|_{2\rightarrow p}> 1$ with $p=3,4$. Note that
\begin{align*}
\|{\widehat{U}}^\top\|_{2\rightarrow p}
\leq& \|{\widehat{U}}^\top\|_{2\rightarrow p}^p=\sup_{\|\x\|_2=1} \|{\widehat{U}}^\top\x\|^p_p.
\end{align*}
Take arbitrary $\x\in \mathbb{S}^{n-1}$,  the unit sphere in $\R^n$. 
Denoting $S$ as the indices of the largest $L$ entries in ${\widehat{U}}^\top\x$, we have
\begin{align*}
\|{\widehat{U}}^\top\x\|_p^p=\|{\widehat{U}}_S^\top\x\|_p^p+\|{\widehat{U}}_{S^c}^\top\x\|_p^p.
\end{align*}

Next, we bound the above two terms sequentially. 

\begin{itemize}
\item    For the first term, we have
\begin{align*}
\|{\widehat{U}}_S^\top\x\|_p^p
&\leq \|{\widehat{U}}_S^\top\x\|_2^2
\\
&\leq \|{\widehat{U}}_S\|^2\cdot \|\x\|_2^2 
\\
&=\|{\widehat{U}}_S{\widehat{U}}_S^\top\| \leq 1+\sum_{i\neq j}|\langle \widehat{\u}_i, \widehat{\u}_j \rangle|
\\
&\leq 1+O\left(\frac{L\gamma}{\sqrt n}\right),
\end{align*}
where the first inequality holds since $({\widehat{U}}_S^\top\x)_i\leq 1$, the fourth inequality follows from the Disk Theorem for symmetric matrix, and the fifth inequality follows from \Cref{asp:incoherence}.

\item
For the second term, note that
\begin{align*}
\min_{i\in S} |\widehat{\u}_i^\top\x|^2
&\leq \frac{1}{L}\|{\widehat{U}}_S{\widehat{U}}^\top_S\| \|\x\|^2
\\
&\leq\frac{1}{L} \bigg(1+O\left(L\frac\gamma{\sqrt n}\right)\bigg)
\\
&= O\left(\frac{1}{L}\right),
\end{align*}
where we assume 
\[\frac{\gamma}{\sqrt n}=o\left(\frac{1}{L}\right)\text{ or }L=o\left(\frac{\sqrt n}{\gamma}\right).
\] 
Choosing a support $S$ such that
$\max_{i\in S^c} |\widehat{\u}_i^\top\x|^2 \leq O(\frac{1}{L}).$
Then, we have
\begin{align*}
\|{\widehat{U}}_{S^c}^\top\x\|_p^p
&= \sum_{i\notin S} |\widehat{\u}_i^\top \x|^p
\\
&\leq \big(\max_{i\notin S} |\widehat{\u}_i^\top\x|^{p-2}\big)\sum_{i\notin S} |\widehat{\u}_i^\top \x|^2 \\
&\leq O(L^{1-0.5p})\|{\widehat{U}}_{S^c}^\top\x\|_2^2
\\
&\leq O(L^{1-0.5p}rn^{-1}),
\end{align*}
where the third inequality holds since
$\max_{i\in S^c} |\widehat{\u}_i^\top\x|^2 \leq O(\frac{1}{L})$, and the last inequality follows from
\[\|{\widehat{U}}_{S^c}^\top\x\|_2^2\leq\|U\|^2\leq(1+c\sqrt{\frac{r}{n}})^2=O(\frac{r}{n}),\]
 which is a consequence of \Cref{asp:spectrum}.
Therefore, when $p=3$, we get
\[\|{\widehat{U}}_{S^c}^\top\x\|_3^3\leq O(L^{-0.5}rn^{-1}).\]

\end{itemize}

Combining these two terms, we have 
$$
\|{\widehat{U}}^\top\x\|_3^3\leq 1+ O(L{\gamma}/{\sqrt n})+O(L^{-0.5}rn^{-1}),
$$
which can be further optimized by carefully choosing $L$ and $r$. 

In particular, let us consider the following optimization:
$$
\min_{r\gg n,L} \max ~ \big\{O(L{\gamma}/{\sqrt n}),O(L^{-0.5}rn^{-1})\big\}.
$$
Choosing
$
L=O({ n^{0.5-r_c}}/\gamma),
$ and $r = O(n^{1+\beta})$,
we get
\begin{align*}
\big\{O(L\gamma/{\sqrt n}),O(L^{-0.5}rn^{-1})\big\}&=\big\{O(n^{-r_c}),O(\gamma^{0.5}n^{-1.25+0.5r_c}r)\big\}\\
&=\big\{O(n^{-r_c}),O(\gamma^{0.5}n^{-0.25+0.5r_c+\beta})\big\}\\
& = \big\{O(n^{-r_c}),O(n^{-r_c}\polylog(n)^{0.5}) \big\},
\end{align*}
where the last equality follows by setting 
$-r_c=-0.25+0.5r_c+\beta \text{ or } \beta=0.25-1.5r_c$.
It follows that
\begin{align*}
\|{\widehat{U}}^\top\x\|_3^3\leq 1+ O(\gamma n^{-r_c}),
\end{align*}
when $r=O( n^{1.25-1.5r_c})$.
Since the above inequality holds for any $\x\in \mathbb{S}^{n-1}$, we have
\begin{align*}
\max_{\x\in \mathbb{S}^{n-1}}\|{\widehat{U}}^\top\x\|_3^3&\leq 1+ O(\operatorname*\gamma n^{-r_c})
\end{align*}
when provided with $r = O( n^{1.25-1.5r_c})$. Then, we obtain
\begin{align*}
\|{\widehat{U}}^\top\|_{2\rightarrow3}^3\leq 1+ O(\operatorname*\gamma n^{-r_c}),    
\end{align*}
which further implies that
\begin{align*}
\|{\widehat{U}}^\top\|_{2\rightarrow3}\leq 1+ O(\operatorname*\gamma n^{-r_c}),    
\end{align*}
when provided with $r = O( n^{1.25-1.5r_c})$.

With a similar argument, we can also show that
\begin{align*}
\|{\widehat{U}}^\top\|_{2\rightarrow4}\leq 1+ O(\polylog(n)n^{-r_c})    
\end{align*}
when provided with $r = O( n^{1.5-2r_c})$.
\end{proof}

\begin{lem}\label{lem:aux:isometry}
Under \Cref{asp:isometry} on $\widehat U$, we have
\begin{align}
\lambda_{\min}(({{U^\star}^\top}{U^\star})\odot({{U^\star}^\top}{U^\star}))
&\geq \underline{c}^4(1-\frac{\gamma\sqrt r}{n}),\nonumber\\
\lambda_{\max}(({{U^\star}^\top}{U^\star})\odot({{U^\star}^\top}{U^\star}))
&\leq \bar{c}^4(1+\frac{\gamma\sqrt r}{n}).
\end{align}
\end{lem}

\begin{proof}
From \Cref{asp:isometry}, we have
$$
\|(\widehat{U}^\top\widehat{U})\odot({\widehat{U}}^\top\widehat{U})-I_r\|\leq\frac{\gamma\sqrt r}{n}.
$$
Since the operator norm of a symmetric matrix equals the largest absolute eigenvalues  and
$$
\max_{i\in[r]}|\lambda_i((\widehat{U}^\top\widehat{U})\odot({\widehat{U}}^\top\widehat{U})-I)|=\max_{i\in[r]}|\lambda_i((\widehat{U}^\top\widehat{U})\odot({\widehat{U}}^\top\widehat{U}))-1|,
$$
 we have
\begin{align*}
\max_{i\in[r]}|\lambda_i((\widehat{U}^\top\widehat{U})\odot({\widehat{U}}^\top\widehat{U}))-1|\leq\frac{\gamma \sqrt r}{n},
\end{align*}
implying
\begin{align*}
|\lambda_i((\widehat{U}^\top\widehat{U})\odot({\widehat{U}}^\top\widehat{U}))-1|\leq\frac{\gamma \sqrt r}{n},i\in[r].
\end{align*}
Hence, we get  
\begin{align*}
1-\frac{\gamma \sqrt r}{n}\leq\lambda_i((\widehat{U}^\top\widehat{U})\odot({\widehat{U}}^\top\widehat{U}))\leq1+\frac{\gamma \sqrt r}{n}.
\end{align*}

Next, we derive the lower bound and upper bound for  $\lambda_i({{U^\star}^\top}U^\star\odot{{U^\star}^\top}U^\star)$. Observe that
\begin{align*}
\bigg[(\widehat{U}^\top\widehat{U})\odot({\widehat{U}}^\top\widehat{U})\bigg]_{i,j}&=(\widehat{\u}_i^\top\widehat{\u}_j)^2,\\
\bigg[({{U^\star}^\top}{U^\star})\odot({{U^\star}^\top}{U^\star})\bigg]_{i,j}&={c_i^\star}^{2/3}{c_j^\star}^{2/3}(\widehat{\u}_i^\top\widehat{\u}_j)^2.
\end{align*}
Then, we have
\begin{align*}
\lambda_{\min}(({{U^\star}^\top}{U^\star})\odot({{U^\star}^\top}{U^\star}))
=&\min_{\|\x\|_2=1} \x^\top(({{U^\star}^\top}{U^\star})\odot({{U^\star}^\top}{U^\star}))\x\\
=&\min_{\|\x\|_2=1} \sum_{i,j\in[r]}{c_i^\star}^{2/3}{c_j^\star}^{2/3} x_ix_j(\widehat{\u}_i^\top\widehat{\u}_j)^2\\
\geq & \min_{i\in[r]} {c_i^\star}^{4/3} \cdot \min_{\|\x\|_2=1} \sum_{i,j\in[r]} x_ix_j(\widehat{\u}_i^\top\widehat{\u}_j)^2\\
=&\underline{c}^4\cdot \lambda_{\min}((\widehat{U}^\top\widehat{U})\odot({\widehat{U}}^\top\widehat{U})).
\end{align*}
Similarly, we can obtain
\begin{align*}
\lambda_{\max}(({{U^\star}^\top}{U^\star})\odot({{U^\star}^\top}{U^\star}))
&\leq\bar{c}^4\cdot \lambda_{\max}((\widehat{U}^\top\widehat{U})\odot({\widehat{U}}^\top\widehat{U})).
\end{align*}
This completes the proof of Lemma~\ref{lem:aux:isometry}.
\end{proof}

\begin{lem}\label{lem:aux:von:Neumann}
For any $n \times n$ symmetric semi-definite matrices $A$, $B$, we have
\begin{align}
\langle A, B\rangle \leq \min\{\|B\| \tmop{trace}(A),\|A\| \tmop{trace}(B)\}.
\end{align}
\end{lem}

\begin{proof}
Due to symmetry, we only prove $\langle A, B\rangle \leq \|B\| \tmop{trace}(A)$.
Denote the eigenvalues  of $A$ as $\lambda_i(A),i\in[n]$ and the eigenvalues of $B$ as $\lambda_i(B),i\in[n]$. Note that for a symmetric semidefinite matrix, the eigenvalues and the singular\
values are the same. Then, by Von Neumann's trace inequality, we obtain
\begin{align*}
\langle A, B\rangle
&\leq \sum_{i\in[n]}\lambda_i(A)\lambda_i(B)
\\
&\leq \lambda_{\max}(B)\sum_{i\in[n]}\lambda_i(A)
\\
&=\|B\|\tmop{trace}(A).
\end{align*}
\end{proof}

\begin{lem}\label{lem:aux:Hardys}
For any $k\geq 2$, we have
$\|AX\|_k\leq\|A\|_{2\rightarrow k}\|X\|_F$.
\end{lem}

\begin{proof}
Observe that
\begin{align*}
\|AX\|_k^k
&\leq\sum_i(\|A\|_{2\rightarrow k}\|\x_i\|_2^k)
\\
&=\|A\|_{2\rightarrow k}^k\sum_i\|\x_i\|_2^k
\\
&\leq\|A\|^k_{2\rightarrow k}\|X\|_F^k,
\end{align*}
where the first inequality follows from the definition of $2\rightarrow k$ norm, and the last inequality holds since
\begin{align*}
\sum_i\|\x_i\|_2^k&=\big\|[\|\x_1\|_2~\cdots~\|\x_r\|_2]\big\|_k^k
\\
&\leq \big\|[\|\x_1\|_2~\cdots~\|\x_r\|_2]\big\|_2^k
\\
&= \|X\|_F^k.
\end{align*}
Here, the second inequality follows from the {\em Hardy's} inequality that $\|\x\|_k\leq\|\x\|_2$ when $k\geq 2$.
Therefore, we have
$\|AX\|_k^k\leq\|A\|^k_{2\rightarrow k}\|X\|_F^k$ or
$\|AX\|_k\leq\|A\|_{2\rightarrow k}\|X\|_F$.
\end{proof}

\begin{lem}\label{lem:aux:3norm}
For any vectors $\mathbf{f},\mathbf{g},\mathbf{h}$ of the same size, we have
\[
\sum_i |\mathbf{f}(i)\mathbf{g}(i)\mathbf{h}(i)|
\leq\|\mathbf{f}\|_3\|\mathbf{g}\|_3\|\mathbf{h}\|_3.
\]
\end{lem}

\begin{proof} 
Denote $\x\doteq\g\odot \h$. We have
\begin{align*}
\sum_i |\mathbf{f}(i)\mathbf{g}(i)\mathbf{h}(i)|
= &|\lg \f, \x \rg| 
\\   
\le &  \|\f\|_3 \|\x\|_{3/2}\\
= & \|\f\|_3  \left( \sum_i |\g(i)|^{3/2}|\h(i)|^{3/2}    \right)^{2/3}
\\
\le & \|\f\|_3   \left( \sqrt{\sum_i |\g(i)|^{3}}\sqrt{ \sum_i |\h(i)|^{3}} \right)^{2/3}
\\
= & \|\f\|_3   \left( \sqrt{\|\g\|_3^3}\sqrt{ \|\h\|_3^3} \right)^{2/3}
\\
=& \|\f\|_3 \|\g\|_3 \|\h\|_3,
\end{align*}
where the both inequalities follow from  the {\em H\"{o}lder's inequality} that $|\lg\a,\b \rg|\le \|\a\|_p\|\b\|_q$ with $(p,q)=(3,3/2)$ and $(p,q)=(2,2)$, respectively.
\end{proof}

Now, we are ready to prove Lemmas~\ref{lem:TTs} and  \ref{lem:TTs_inner}.

\subsection{Proof of Lemma~\ref{lem:TTs}}
\label{sec:proof_lem_TTs}

\subsubsection{Upper bound for $\|T-T^\star\|_F^2$}

First of all, we expand 
\[T=U\circ U\circ U=(U^\star+H)\circ(U^\star+H)\circ(U^\star+H)\]
into 8 terms as
\begin{align*}
T-T^\star
&=(U^\star+H)\circ(U^\star+H)\circ(U^\star+H)-U^\star\circ U^\star\circ U^\star\\
&= H\circ U^\star \circ U^\star + U^\star \circ H\circ U^\star + U^\star \circ U^\star \circ H +  H\circ H\circ U^\star + H\circ U^\star \circ H + U^\star \circ H\circ H +  H\circ H\circ H.
\end{align*}
We then plug it into $\|T-T^\star\|_F^2$ and continue to  expand it. With some simplification, we have
\begin{align*}
\|T-T^\star\|^2_F 
=&\|\underbrace{\big(H\circ U^\star\circ U^\star+\cdots\big)}_{\#3} 
 +\underbrace{\big(H\circ H\circ U^\star+\cdots\big)}_{\#3}
 +H\circ H\circ H\|_F^2\\
=&3\langle H^\top H,{U^\star}^\top  U^\star\odot{U^\star}^\top  U^\star\rangle
\tag{Term-(1)}\\
& +6\langle {U^\star}^\top U^\star,(H^\top U^\star)\odot({U^\star}^\top H)
\tag{Term-(2)}\\
 &+6\langle H^\top {U^\star}, ({U^\star}^\top H)^{\odot2}\rangle
 \tag{Term-(3)}\\
 &+12\langle {U^\star}^\top H, H^\top H\odot {U^\star}^\top U^\star \rangle
\tag{Term-(4)} \\
 &+6\langle H^\top H, ({U^\star}^\top H)^{\odot2}\rangle
 \tag{Term-(5)}\\
 &+3\langle{U^\star}^\top U^\star,(H^\top H)\odot(H^\top H) \rangle
 \tag{Term-(6)}\\
 &+6\langle H^\top H,(H^\top U^\star)\odot({U^\star}^\top H)\rangle
\tag{Term-(7)} \\
 &+6\langle {U^\star}^\top H,(H^\top H)\odot(H^\top H)\rangle
 \tag{Term-(8)}\\
 &+\langle H^\top H,(H^\top H)\odot(H^\top H)\rangle,
 \tag{Term-(9)}
\end{align*}
where we use $U^{\odot2}$ to denote $U\odot U$.
Next, we bound the above nine terms sequentially. 

\begin{itemize}[leftmargin=*]
\item \noindent{\em Term-(1):}
With Lemmas \ref{lem:aux:isometry} and \ref{lem:aux:von:Neumann}, the first term can be bounded with
\begin{align*}
\langle H^\top H,{U^\star}^\top  U^\star\odot{U^\star}^\top  U^\star\rangle
&\leq \lambda_{\max}({U^\star}^\top  U^\star\odot{U^\star}^\top  U^\star)\|H\|_F^2
\\
&\leq \bar{c}^4\bigg(1+\frac{\gamma \sqrt r}{n}\bigg)\|H\|_F^2.
\end{align*}

\item     \noindent{\em Term-(2):}
Note that
\begin{align*}
\langle {U^\star}^\top U^\star,(H^\top U^\star)\odot({U^\star}^\top H)
=&\sum_{i,j\in[r]}
\langle\u^\star_i,\u^\star_j\rangle \langle \h_i,\u^\star_j\rangle
\langle \u^\star_i,\h_j\rangle\\
=&\underbrace{\sum_{i\neq j}\langle\u^\star_i,\u^\star_j\rangle \langle \h_i,\u^\star_j\rangle\langle \u^\star_i,\h_j\rangle}_{\Pi_1}
+\underbrace{\sum_{i\in[r]} \langle\u^\star_i,\u^\star_i\rangle \langle \h_i,\u^\star_i\rangle
 \langle \h_i,\u^\star_i\rangle}_{\Pi_2}.
\end{align*}
We first bound $\Pi_1$ with
\begin{align*}
\Pi_1
&\leq\sum_{i\neq j}\big|\langle\u^\star_i,\u^\star_j\rangle \langle \h_i,\u^\star_j\rangle\langle \u^\star_i,\h_j\rangle\big|
\\
&\leq\max_{i\neq j}\big|\langle\u^\star_i,\u^\star_j\rangle\big|\cdot\sum_{i\neq j}\big|\langle \h_i,\u^\star_j\rangle\langle \u^\star_i,\h_j\rangle\big|\\
&\leq\mu^\star\cdot\langle\big|H^\top U^\star\big|,\big|{U^\star}^\top H\big|\rangle
\\
&\leq \mu^\star\cdot\|H^\top U^\star\|_F\|{U^\star}^\top H\|_F
\\
&\leq\mu^\star\cdot\|U^\star\|^2\cdot\|H\|_F^2,
\end{align*}
where the operation $|\cdot|$ on a matrix means taking the absolute value of all its entries.
 The last second inequality follows from the {\em Cauchy-Schwarz's inequality}.
The last inequality holds due to Lemma \ref{lem:aux:Hardys}.
We then bound $\Pi_2$ with
\begin{align*}
\Pi_2
&\leq\max_{i\in[r]} \langle\u^\star_i,\u^\star_i\rangle\cdot\sum_{i\in[r]}\langle \h_i,\u^\star_i\rangle^2
\\
&\leq\max_{i\in[r]} \|\u^\star_i\|_2^2\cdot\sum_{i\in[r]}\|\h_i\|_2^2\|\u^\star_i\|_2^2 \\
&\leq\max_{i\in[r]} \|\u^\star_i\|_2^4\cdot\sum_{i\in[r]}\|\h_i\|_2^2
\\
&=\bar{c}^4\cdot\|H\|_F^2.
\end{align*}
Combining the upper bound for $\Pi_1$ and $\Pi_2$, we obtain
\begin{align*}
\langle {U^\star}^\top U^\star,(H^\top U^\star)\odot({U^\star}^\top H)\leq(\mu^\star\|U^\star\|^2+\bar{c}^4)\cdot\|H\|_F^2.
\end{align*}

\item \noindent{\em Term-(3):}
With Lemmas \ref{lem:aux:Hardys} and \ref{lem:aux:3norm}, we can bound the third term with
\begin{align*}
\langle H^\top {U^\star}, ({U^\star}^\top H)\odot({U^\star}^\top H)\rangle
&\leq \|{U^\star}^\top H\|_3^3
\\
&\leq \|{U^\star}^\top \|_{2\rightarrow 3}^3 \|H\|_F^3.
\end{align*}

\item \noindent{\em Term-(4):}
Note that
\begin{align*}
\langle {U^\star}^\top H, H^\top H\odot {U^\star}^\top U^\star \rangle
&=\sum_{i\neq j} \langle \u^\star_i,\h_j \rangle\langle\h_i,\h_j \rangle\langle\u^\star_i,\u^\star_j \rangle
+\sum_{i} \langle \u^\star_i,\h_i \rangle\langle\h_i,\h_i \rangle\langle\u^\star_i,\u^\star_i \rangle
\\
&\doteq \Pi_1+\Pi_2.
\end{align*}
With a similar technique used in \noindent{\em Term-(2)}, we can bound $\Pi_1$ and $\Pi_2$ with
\begin{align*}
\Pi_1
&\leq\sum_{i\neq j}\big|\langle \u^\star_i,\h_j \rangle\langle\h_i,\h_j \rangle\langle\u^\star_i,\u^\star_j \rangle\big|
\\
&\leq\max_{i\neq j}\big|\langle\u^\star_i,\u^\star_j\rangle\big|\cdot\sum_{i\neq j}\big|\langle \h_i,\h_j\rangle\langle \u^\star_i,\h_j\rangle\big|\\
&=\mu^\star\cdot\langle\big|H^\top H\big|,\big|{U^\star}^\top H\big|\rangle
\\
&\leq\mu^\star\cdot\|U^\star\|\|H\|_F^3,
\end{align*}
and
\begin{align*}
\Pi_2
&\leq\max_{i\in[r]} \langle\u^\star_i,\u^\star_i\rangle\cdot\sum_{i\in[r]}\langle \u^\star_i,\h_i \rangle\langle\h_i,\h_i \rangle
\\
&\leq\max_{i\in[r]} \|\u_i\|_2^2\cdot\sum_{i\in[r]}\|\h_i\|_2^3\|\u^\star_i\|_2 \\
&\leq\max_{i\in[r]} \|\u_i\|_2^3\cdot\sum_{i\in[r]}\|\h_i\|_2^3
\\
&\leq\bar{c}^3\cdot\|H\|_F^3,
\end{align*}
where the second inequality follows from the {\em Cauchy-Schwarz's inequality}, and the last inequality follows from the {\em Hardy's inequality}, i.e.,
\begin{align*}
\sum_{i\in[r]}\|\h_i\|_2^3
=&\big\|[\|\h_1\|_2~\cdots~\|\h_r\|_2]\big\|_3^3\leq \big\|[\|\h_1\|_2~\cdots~\|\h_r\|_2]\big\|_2^3=\|H\|_F^3
\end{align*}
Then, we have
\begin{align*}
\langle {U^\star}^\top H, H^\top H\odot {U^\star}^\top U^\star \rangle &\leq (\mu^\star\|U^\star\|+\bar{c}^3)\|H\|_F^3.
\end{align*}

\item \noindent{\em Term-(5):}
The fifth term can be bounded with
\begin{align*}
\langle H^\top H, ({U^\star}^\top H)^{\odot2}\rangle
\leq& \|H^\top H\|_F\|({U^\star}^\top H)^{\odot2}\|_F\\
=&\|H^\top H\|_F\bigg(\|{U^\star}^\top H\|_4^4\bigg)^{1/2}
\\
\leq&\|{U^\star}^\top \|_{2\rightarrow4}^2\|H\|_F^4,
\end{align*}
where the first inequality results from the {\em Cauchy-Schwarz's inequality}, and the last inequality follows from Lemma \ref{lem:aux:Hardys} and the fact that $\|H\|\leq\|H\|_F$ holds for any matrix $H$.

\item \noindent{\em Term-(6)}:
The sixth term can be rewritten as
\begin{align*}
\langle{U^\star}^\top U^\star,(H^\top H)\odot(H^\top H) \rangle
=&\sum_{i,j} \langle \h_i,\h_j \rangle
\langle\u_i,\u_j \rangle\langle\h_i,\h_j
\rangle
\\
=&\sum_{i\neq j} \langle \h_i,\h_j \rangle \langle\u^\star_i,\u^\star_j \rangle\langle\h_i,\h_j \rangle
+\sum_{i} \langle \h^\star_i,\h_i \rangle \langle\u^\star_i,\u^\star_i \rangle\langle\h_i,\h_i \rangle
\\
\doteq& \Pi_1+\Pi_2.
\end{align*}
Note that we can bound $\Pi_1$ as
\begin{align*}
\Pi_1
&\leq\sum_{i\neq j}\big|\langle \h_i,\h_j \rangle\langle\u^\star_i,\u^\star_j \rangle\langle\h_i,\h_j \rangle\big|
\\
&\leq\max_{i\neq j}\big|\langle\u^\star_i,\u^\star_j\rangle\big|\cdot\sum_{i\neq j}\langle \h_i,\h_j\rangle^2 \\
&\leq\mu^\star\cdot\langle\big|H^\top H\big|,\big|H^\top H\big|\rangle
\\
&\leq\mu^\star\|H\|_F^4,
\end{align*}
where the third inequality results from the {\em Cauchy-Schwarz's inequality}, and the fifth inequality follows from Lemma~\ref{lem:aux:Hardys}.
We can also bound $\Pi_2$ as
\begin{align*}
\Pi_2
&\leq\max_{i\in[r]} \langle\u^\star_i,\u^\star_i\rangle\cdot\sum_{i\in[r]}\langle \h_i,\h_i \rangle \langle\h_i,\h_i \rangle
\\
&\leq\max_{i\in[r]} \|\u_i\|_2^2\cdot\sum_{i\in[r]}\|\h_i\|_2^4
\\
&\leq\bar{c}^2\cdot\|H\|_F^4,
\end{align*}
where the last inequality is due to the {\em Hardy's inequality}, i.e.,
\begin{align*}
\sum_{i\in[r]}\|\h_i\|_2^4
=&\big\|[\|\h_1\|_2~\cdots~\|\h_r\|_2]\big\|_4^4
\\
\leq& \big\|[\|\h_1\|_2~\cdots~\|\h_r\|_2]\big\|_2^4
\\
= &\|H\|_F^4.
\end{align*}
Combining the bound of $\Pi_1$ and $\Pi_2$, we get
\begin{align*}
\langle{U^\star}^\top U^\star,(H^\top H)\odot(H^\top H) \rangle
\leq (\mu^\star+\bar{c}^2)\cdot\|H\|_F^4.
\end{align*}

\item \noindent{\em Term-(7):}
Note that
\begin{align*}
\langle H^\top H,(H^\top U^\star)\odot({U^\star}^\top H)\rangle
\leq& \|H^\top H\|_3\|H^\top U^\star\|_3\|{U^\star}^\top H\|_3\\
\leq &\|H^\top H\|_F\|{U^\star}^\top H\|_3^2
\\
\leq &\|{U^\star}^\top \|_{2\rightarrow 3}^2\|H\|_F^4,
\end{align*}
where the first inequality follows from Lemma \ref{lem:aux:3norm}, the second inequality follows from the {\em Hardy's inequality}, i.e., $|H^\top H\|_3\leq|H^\top H\|_F$, and the last inequality follows from Lemma \ref{lem:aux:Hardys} and $\|H\|\leq\|H\|_F$.

\item \noindent{\em Term-(8):}
With Lemmas Lemma \ref{lem:aux:Hardys} and Lemma \ref{lem:aux:3norm}, we can bound the eighth term as
\begin{align*}
\langle {U^\star}^\top H,(H^\top H)\odot(H^\top H)\rangle
&\leq\|{U^\star}^\top H\|_3\|H^\top H\|_3\|H^\top H\|_3
\\
&\leq\|{U^\star}^\top \|_{2\rightarrow3}\|H\|_F^5,
\end{align*}
where the second inequality follows from the {\em Hardy's inequality}.

\item \noindent{\em Term-(9):}
Similarly, we can bound the last term with
\begin{align*}
\langle H^\top H,(H^\top H)\odot(H^\top H)\rangle
&\leq\|H^\top H\|_3^3\leq\|H^\top H\|_F^3\leq\|H\|_F^6.
\end{align*}
\end{itemize}
 
 \bigskip

\paragraph{Putting together}
By combining the above upper bound for the nine terms, we can obtain
\begin{align*}
\|T-T^\star\|_F^2
\leq &\bigg( 3\bar{c}^4\bigg(1+\frac{\gamma \sqrt r}{n}\bigg)+6(\mu^\star\|U^\star\|^2+\bar{c}^4)\bigg)\|H\|_F^2
\\
&+\bigg(6\|{U^\star}^\top \|_{2\rightarrow3}^3+12(\mu^\star\|U^\star\|+\bar{c}^4  )\bigg)\|H\|_F^3\\
&+\bigg(6\|{U^\star}^\top \|_{2\rightarrow4}^2+6\|{U^\star}^\top \|_{2\rightarrow3}^2+3(\mu^\star+\bar{c}^2)\bigg)\|H\|_F^4
\\
&+6\|{U^\star}^\top \|_{2\rightarrow3}\|H\|_F^5+\|H\|_F^6.
\end{align*}

Note that the objective of this paper is to study the asymptotic performance of the gradient descent algorithm when applied to the tensor decomposition problem, which means that we consider the case when $n$ is sufficiently large. Thus, we will bound $\|T-T^\star\|_F^2$ in the case when $n\rightarrow\infty$. Observe that the upper bound of $\|T-T^\star\|_F^2$ involves the terms $\mu^\star\|U^\star\|$, $\mu^\star\|U^\star\|^2$, $\|{U^\star}^\top \|_{2\rightarrow 3}^k$ and $\|{U^\star}^\top \|_{2\rightarrow 4}^k$.
To bound $\|T-T^\star\|_F^2$ in the asymptotic sense, we next compute the asymptotic upper bound of these involved terms.

Using Lemmas \ref{lem:aux:incoherence} and \ref{lem:aux:spectrum}, we get
\begin{align*}
\mu^\star\|U^\star\|&\leq \frac{\bar{c}^3\gamma}{\sqrt n}\bigg(1+c_1\sqrt{\frac{r}{n}}\bigg)\overset{n\rightarrow \infty}{=}0,
\end{align*}
where the last equality holds provided that $r\ll n^2$. 

Similarly, we can also get
\begin{align*}
\mu^\star\|U^\star\|^2&\leq \frac{\bar{c}^4\gamma}{\sqrt n}\bigg(1+c_1\sqrt{\frac{r}{n}}\bigg)^2\overset{n\rightarrow \infty}{=}0,
\end{align*}
when provided with $r\ll n^{1.5}$.
Recall that  
\begin{align*}
&\|{U^\star}^\top \|^k_{2\rightarrow 3}
\leq \bar{c}^k\big(1+O(\gamma n^{-\epsilon})\big), \quad \tmop{if}~r=O(n^{1.25-1.5\epsilon}),
\\
&\|{U^\star}^\top \|^k_{2\rightarrow 4}
\leq \bar{c}^k\big(1+O(\gamma n^{-\epsilon})\big), \quad \tmop{if}~r=O( n^{1.5-2\epsilon}), 
\end{align*}
which implies that when $n$ goes to infinity:
\begin{align*}
&\|{U^\star}^\top \|^k_{2\rightarrow 3}
\leq \bar{c}^k, \quad \tmop{if}~r=O(n^{1.25-1.5\epsilon}),
\\
&\|{U^\star}^\top \|^k_{2\rightarrow 4}
\leq \bar{c}^k, \quad \tmop{if}~r=O( n^{1.5-2\epsilon}). 
\end{align*}

Combining these asymptotic bound and letting $r=O(n^{1.25-1.5\epsilon})$, we then have 
\begin{align*}
\|T-T^\star\|_F^2
&\leq
\bigg(9\bar{c}^4+18\bar{c}^3\|H\|_F+15\bar{c}^2\|H\|_F^2
+6\bar{c}\|H\|_F^3+\|H\|_F^4
\bigg)\|H\|_F^2
\\
&\leq
\bigg(9\omega^4\underline{c}^4+18\omega^3\underline{c}^3\|H\|_F
+15\omega^2\underline{c}^2\|H\|_F^2
+6\omega\underline{c}\|H\|_F^3+\|H\|_F^4
\bigg)\|H\|_F^2
\end{align*}
holds for sufficiently large $n$. Here, the last inequality follows from the fact that $\omega\geq 1$.

Finally, setting \[\|H\|_F\leq0.07\frac{\underline{c}}{\omega^3},\]
 we get
\begin{align*}
\|T-T^\star\|_F^2&\leq \bigg(9+\frac{1}{\omega^4}\bigg(18(0.07)+15\frac{(0.07)^2}{\omega^4}+6\frac{(0.07)^3}{\omega^8}
+\frac{(0.07)^4}{\omega^{12}}
\bigg)\bigg)\omega^4\underline{c}^4\|H\|_F^2 \\
&\leq 10.336\omega^4\underline{c}^4\|H\|_F^2.
\end{align*}
In conclusion, under  Assumptions \ref{asp:incoherence}-\ref{asp:isometry}  if $r=O(n^{1.25})$ and $\|H\|_F\leq0.07\frac{\underline{c}}{\omega^3}$, we have
$$
\|T-T^\star\|_F^2\leq
 10.336\omega^4\underline{c}^4\|H\|_F^2
$$
holds for sufficiently large $n$.

\subsubsection{Lower bound for  $\|T-T^\star\|_F^2$}

Similar to computing upper bound of the nine terms, we next compute the lower bound of these nine terms sequentially. 

\begin{itemize}[leftmargin=*]
\item    \noindent{\em Term-(1):}
The first term can be bounded with
\begin{align*}
\langle H^\top H,{U^\star}^\top  U^\star\odot{U^\star}^\top  U^\star\rangle
&\geq \lambda_{\min}({U^\star}^\top  U^\star\odot{U^\star}^\top  U^\star)\|H\|_F^2.
\end{align*}

\item \noindent{\em Term-(2):}
Observe that
\begin{align*}
\langle {U^\star}^\top U^\star,(H^\top U^\star)\odot({U^\star}^\top H)
&=\sum_{i,j\in[r]}
\langle\u^\star_i,\u^\star_j\rangle \langle \h_i,\u^\star_j\rangle
\langle \u^\star_i,\h_j\rangle\\
&\geq \sum_{i\neq j}
\langle\u^\star_i,\u^\star_j\rangle \langle \h_i,\u^\star_j\rangle
\langle \u^\star_i,\h_j\rangle\\
&\geq  -\max_{i\neq j}|\langle\u^\star_i,\u^\star_j\rangle|\cdot\sum_{i\neq j}|\langle\h_i,\u^\star_j\rangle|
\cdot|\langle\u^\star_i,\h_j\rangle|\\
&\geq -\mu^\star\bigg(\sum_{i,j\in[r]}\langle\h_i,\u^\star_j\rangle^2\bigg)^{1/2}
\bigg(\sum_{i,j}\langle\h_j,\u^\star_i\rangle^2\bigg)^{1/2}\\
&= - \mu^{\star} \|H^\top U^\star\|_F \cdot\|{U^\star}^\top H\|_F\\
&\geq   -\mu^{\star}  \|U^\star\|^2\cdot\|H\|_F^2,
\end{align*}
where the second inequality follows from
$\sum_{i\in[r]} \langle\u^\star_i,\u^\star_i\rangle \langle \h_i,\u^\star_i\rangle\langle \u^\star_i,\h_i\rangle\geq 0$, the fourth inequality follows from the {\em Cauchy-Schwarz's inequality}, and the last inequality is a consequence of Lemma~\ref{lem:aux:Hardys}.

\item \noindent{\em Term-(3):}
With Lemmas~\ref{lem:aux:Hardys} and \ref{lem:aux:3norm}, we can bound the third term as
\begin{align*}
\langle H^\top {U^\star}, ({U^\star}^\top H)^{\odot2}\rangle
&\geq -\|{U^\star}^\top H\|_3^3
\\
&\geq -\|{U^\star}^\top \|_{2\rightarrow 3}^3 \|H\|_F^3.
\end{align*}

\item \noindent{\em Term-(4):}
With a similar technique used in {\em Term-(2)}, the forth term can be bounded with
\begin{align*}
\langle {U^\star}^\top H, H^\top H\odot {U^\star}^\top U^\star \rangle
&=\sum_{i\neq j} \langle \u^\star_i,\h_j \rangle
\langle\h_i,\h_j \rangle\langle\u^\star_i,\u^\star_j \rangle
+\sum_{i} \langle \u^\star_i,\h_j \rangle
\langle\h_i,\h_i \rangle\langle\u^\star_i,\u^\star_i \rangle\\
&\geq -\mu^\star\sum_{i\neq j}|\langle \u^\star_i,\h_j \rangle\langle\h_i,\h_j \rangle|-\bar{c}^3\bigg(\sum_i\|\h_i\|_2^3\bigg)\\
&\geq -\mu^\star\|U^\star\|\|H\|_F^3-\bar{c}^3\|H\|_F^3\\
&=-(\mu^\star\|U^\star\|+\bar{c}^3)\|H\|_F^3.
\end{align*}

\item \noindent{\em Term-(5):}
Using {\em Cauchy-Schwarz's inequality} and Lemma~\ref{lem:aux:Hardys}, we have
\begin{align*}
\langle H^\top H, ({U^\star}^\top H)^{\odot2}\rangle
\geq& -\|H^\top H\|_F\|({U^\star}^\top H)^{\odot2}\|_F\\
=&-\|H^\top H\|_F\bigg(\|{U^\star}^\top H\|_4^4\bigg)^{1/2}
\\
\geq &-\|{U^\star}^\top \|_{2\rightarrow4}^2\|H\|_F^4.
\end{align*}

\item \noindent{\em Term-(6)}:
The sixth term can be simply bounded with
\begin{align*}
\langle{U^\star}^\top U^\star,(H^\top H)\odot(H^\top H) \rangle\geq0.
\end{align*}

\item
{\em Term-(7):}
Observe that
\begin{align*}
\langle H^\top H,(H^\top U^\star)\odot(U^{\star'}H)\rangle
&\geq -\|H^\top H\|_3\|H^\top U^\star\|_3\|{U^\star}^\top H\|_3\\
&\geq -\|H^\top H\|_F\|{U^\star}^\top H\|_3^2\\
&\geq -\|{U^\star}^\top \|_{2\rightarrow 3}^2\|H\|_F^4,
\end{align*}
where the first line follows from Lemma \ref{lem:aux:3norm}, the second line holds due to the {\em Hardy's inequality}, i.e., $\|H^\top H\|_F\geq\|H^\top H\|_3$, and the last inequality is a consequence of Lemma~\ref{lem:aux:Hardys}.

\item \noindent{\em Term-(8):}
Note that
\begin{align*}
\langle {U^\star}^\top H,(H^\top H)\odot(H^\top H)\rangle
&\geq-\|{U^\star}^\top H\|_3\|H^\top H\|_3\|H^\top H\|_3 \\
&\geq-\|{U^\star}^\top H\|_3\|H^\top H\|_F\|H^\top H\|_F \\
&\geq-\|{U^\star}^\top \|_{2\rightarrow3}\|H\|_F^5.
\end{align*}
{\em Term-(9):}
The last term can be simply bounded with
\begin{align*}
\langle H^\top H,(H^\top H)\odot(H^\top H)\rangle\geq0.
\end{align*}

\end{itemize}

\bigskip

\paragraph{Putting together}
Combining the above lower bound for the nine terms, we can obtain
\begin{align*}
\|T-T^\star\|_F^2
\geq &
\bigg(3\lambda_{\min}({U^\star}^\top  U^\star\odot{U^\star}^\top  U^\star)-6\mu^\star\|U^\star\|^2\bigg)\|H\|_F^2
\\
&-\bigg(6\|{U^\star}^\top \|_{2\rightarrow 3}^3+12(\mu^\star\|U^\star\|+\bar{c}^3)\bigg)\|H\|_F^3\\
&-\bigg(6\|{U^\star}^\top \|_{2\rightarrow 4}^2+6\|{U^\star}^\top \|_{2\rightarrow 3}^2\bigg)\|H\|_F^4
\\
&-\bigg(6\|{U^\star}^\top \|_{2\rightarrow 3}\bigg)\|H\|_F^5.
\end{align*}
Observe that the above lower-bound involves terms $\mu^\star\|U^\star\|$, $\mu^\star\|U^\star\|^2$, $\|U^\star\|_{2\rightarrow 3}^k$,
$\|U^\star\|_{2\rightarrow 4}^k$ and $\lambda_{\min}({U^\star}^\top  U^\star\odot{U^\star}^\top  U^\star)$.
Since we have studied the asymptotic bound of the first four terms, it remains to give the asymptotic bound of the last term
$\lambda_{\min}({U^\star}^\top  U^\star\odot{U^\star}^\top  U^\star)$. With Lemma \ref{lem:aux:isometry}, we have
$$\lambda_{\min}(({U^\star}^\top {U^\star})\odot({U^\star}^\top {U^\star}))\geq \underline{c}^4(1-\frac{\gamma\sqrt r}{n}).$$
Then, as $n \rightarrow \infty$, we have
$$\lambda_{\min}(({U^\star}^\top {U^\star})\odot({U^\star}^\top {U^\star}))\geq \underline{c}^4.$$
Plugging all the asymptotic bound of the five terms, we obtain
\begin{align*}
\|T-T^\star\|_F^2
&\geq(3\underline{c}^4-18\bar{c}^3\|H\|_F-12\bar{c}^2\|H\|_F^2-6\bar{c}\|H\|_F^3)\|H\|_F^2\\
&\geq(3\underline{c}^4-18\omega^3\underline{c}^3\|H\|_F-12\omega^2\underline{c}^2\|H\|_F^2-6\omega\underline{c}\|H\|_F^3)\|H\|_F^2,
\end{align*}
where the second line follows from $\omega>1$.

Finally, setting 
\[\|H\|_F\leq0.07\frac{\underline{c}}{\omega^3},\] we  get
$$
\|T-T^\star\|_F^2
\geq\bigg( 3 -\bigg(18(0.07)+\frac{12(0.07)^2}{\omega^4}+\frac{6(0.07)^3}{\omega^{8}}\bigg) \bigg)\underline{c}^4\|H\|_F^2
\geq 1.679\underline{c}^4\|H\|_F^2.
$$

\section{Proof of Lemma~\ref{lem:TTs_inner}}
\label{sec:proof_lem_TTs_inner}

We need the following lemmas to prove Lemma~\ref{lem:TTs_inner}.

\begin{lem}\label{lem4:1}
For any matrix $U$, we have
$
\|U\circledast U\|\leq \|U\|^2.
$ 

\end{lem}
\begin{proof}
Note that
\begin{align*}
\|U\circledast U\|
&= \max_{\|\bxi\|_2=1} \|\sum_l \xi_l \u_l\circledast \u_l \|_2
\\
&=\max_{\|\bxi\|_2=1} \|\sum_l \xi_l \u_l \u_l^\top\|_F\\
&=\max_{\|\bxi\|_2=1} \|U\tmop{diag}(\boldsymbol{\bxi})U^\top \|_F
\\
&\leq\|U\|^2\cdot\max_{\|\bxi\|_2=1}\|\tmop{diag}(\boldsymbol{\bxi})\|_F
\\
&=\|U\|^2,
\end{align*}
where the first inequality follows from Lemma \ref{lem:aux:Hardys}, and the last equality holds since $\|\tmop{diag}(\boldsymbol{\bxi})\|_F=\|\bxi\|_2=1$.
\end{proof}

\begin{lem}\label{lem4:2}
For any  symmetric tensor $Q\in\R^{n\times n\times n}$,  we have
$
\|Q\|\leq\|Q_{(1)}\|.
$
\end{lem}
\begin{proof}
This is due to their different degrees of freedom in optimizing variable. To be more precise, with the definition of tensor operator norm, we have
\begin{align*}
\|Q\|&=\max_{\bxi\in\S^{n-1}} Q \bar{\times}_1 \bxi \bar{\times}_2 \bxi \bar{\times}_3 \bxi
\\
&=\max_{\bxi\in\S^{n-1}}  \bxi^\top Q_{(1)} (\bxi\circledast\bxi)\\
&\le\max_{\bxi\in\S^{n-1},\y\in\S^{n^2-1}} \bxi^\top Q_{(1)}\y\\
&=\|Q_{(1)}\|,
\end{align*}
where the inequality follows from the feasible set of the optimization problem defining $\|Q_{(1)}\|$ covers that of the optimization problem defining $\|Q\|$, and the objective functions of the both optimization problems are of the same form.
\end{proof}

\begin{lem}\label{lem4:3}
For a differentiable scalar function $f$ defined on symmetric tensors, we have
\begin{align*}
\|\nabla_U f(U\circ U\circ U)\|\leq 3\|[\nabla f(U\circ U\circ U)]_{(1)} \|\cdot \|U\|^2,
\end{align*}
or simply write
$\|\nabla_U f\|\leq 3\|[\nabla f]_{(1)}\|\cdot \|U\|^2$.
\end{lem}

\begin{proof}
Note that
\begin{align*}
\|\nabla_U f(U\circ U\circ U)\|&=\|3[\nabla f(U\circ U\circ U)]_{(1)}  (U\circledast U)\| \\
&\leq 3\|[\nabla f(U\circ U\circ U)]_{(1)} \|\cdot\|U\circledast U\|\\
&\leq 3 \|[\nabla f(U\circ U\circ U)]_{(1)} \|\cdot \|U\|^2,
\end{align*}
where the first equality comes from $\nabla_U f(U\circ U\circ U)=3[\nabla f(U\circ U\circ U)]_{(1)}  (U\circledast U)$, the first inequality follows from the definition of operator norm, and the last inequality from Lemma \ref{lem4:1}.
\end{proof}

\begin{lem}\label{lem4:4}
For a differentiable scalar function $f$ defined on symmetric tensors. Assume $U=[\u_1~\cdots~\u_r]$.  Then we have
\begin{align*}
\max_l\|(\nabla_Uf)_l\|_2\leq 3\|\nabla f\|\cdot \max_l\|\u_l\|_2^2.
\end{align*}
\end{lem}
\begin{proof}
Direct computation gives
$$
(\nabla_Uf)_l=3\nabla f \times_1 \u_l \times_2 \u_l.
$$
Then,
\begin{align*}
\max_l\|(\nabla_Uf)_l\|_2
&=3\max_l\nabla f \times_1 \u_l \times_2 \u_l
\leq\|\nabla f\|\cdot \max_l\|\u_l\|_2^2.
\end{align*}
The last  follows from the definition of tensor operator norm.
\end{proof}

\begin{proof}[Proof of Lemma \ref{lem:TTs_inner}] We denote $\triangle = \nabla_U f(U\circ U\circ U)$ through all the proof.
Define $T^+=U^+ \circ U^+ \circ U^+$ and assume that $T$ is an arbitrary symmetric tensor with rank $r$.
Since $f$ is $(r,m,M)$-restricted strongly convex and smooth, we have
\begin{align*}
f(T) &\geq f(T^+) - \langle \nabla f(T), T^+-T \rangle - \frac{M}{2}\|T^+-T\|_F^2  \\
&\geq  f(T^\star) - \langle \nabla f(T), T^+-T \rangle - \frac{M}{2}\|T^+-T\|_F^2, \\
f(T^\star) &\geq f(T) + \langle \nabla f(T), T^\star-T \rangle + \frac{m}{2}\|T^\star-T\|_F^2,
\end{align*}
where the second line holds since $f(T^\star)$ is the global minimum of $f$ hence $f(T^\star)\leq f(T^+)$.
Combining the above two inequalities, we get
\begin{align}\label{lem4:final:1}
\langle \nabla f(T), T-T^\star \rangle\geq \langle \nabla f(T), T-T^+\rangle -\frac{M}{2}\|T-T^+\|_F^2+\frac{m}{2}\|T^\star-T\|_F^2.
\end{align}
Therefore, to bound $\langle \nabla f(T), T-T^\star \rangle$, it suffices to compute a lower bound of $\langle \nabla f(T), T-T^\star \rangle$ and an upper bound of $M\|T-T^+\|_F^2$.

Intuitively, as long as the step size is sufficiently small, $\langle \nabla f(T), T-T^\star \rangle$ would be well-bounded.
Hence, it is critical to choose a stepsize to bound $\langle \nabla f(T), T-T^\star \rangle$. We use the following three rules to
choose the stepsize:
\begin{align*}
 \textbf{Rule I}~~~&3\eta \|\nabla f\| \max_{l\in[r]} \|\u_l\|_2\leq\frac{1}{L_1}, \\
  \textbf{Rule II}~~~&\eta \|U\|^2 \max_{l\in[r]} \|\u_l\|_2^2\leq\frac{1}{ML^2_2},\\
  \textbf{Rule III}~~~&3\eta \|[\nabla f]_{(1)}\|\cdot\|U\|\leq\frac{1}{L_3},
\end{align*}   
where $U$ is the current iteration variable, $M$ is the smoothness coefficient of $f$ and the constant $L_1, L_2, L_3$ will be carefully determined later.

Next, we bound $\langle \nabla f(T), T-T^\star \rangle$ and $M\|T-T^+\|_F^2$ in sequence. Note that
$$T^+=(U-\eta\triangle)\circ(U-\eta\triangle)\circ(U-\eta\triangle),$$ and
\begin{align*}
T-T^+=&\eta(U\circ U\circ\triangle+U\circ\triangle\circ U+\triangle\circ U\circ U)
-\eta^2(U\circ\triangle\circ\triangle+\triangle\circ U\circ \triangle+\triangle\circ \triangle\circ U)+\eta^3 \triangle\circ\triangle\circ\triangle,
\end{align*}
which implies that $\langle \nabla f(T), T-T^+\rangle$ includes seven terms. Since $\nabla f(T)$ is a symmetric tensor, these seven terms can be classified into three classes:
\begin{align*}
&\textbf{(1) }~~~\langle\nabla f(T), U\circ U\circ \triangle\rangle
=\langle \nabla f(T), U\circ\triangle\circ U\rangle=\langle \nabla f(T), \triangle\circ U\circ U\rangle,\\
&\textbf{(2) }~~~\langle\nabla f(T), U\circ\triangle\circ\triangle\rangle
=\langle \nabla f(T), \triangle\circ U\circ\triangle\rangle=\langle \nabla f(T),  \triangle\circ\triangle\circ U\rangle,\\
&\textbf{(3) }~~~\langle\nabla f(T), \triangle\circ \triangle\circ \triangle\rangle.
\end{align*}

For the first class \textbf{(1)}, note that
\begin{align*}
\langle \nabla f(T), \eta(U\circ U\circ \triangle+\cdots)\rangle
=&3\eta \langle\nabla f, U\circ U\circ\triangle\rangle \\
=&\eta \langle 3[\nabla f]_{(1)}~U\circledast U, \triangle\rangle
\\
=&\eta \langle \triangle, \triangle\rangle
\\
=&\eta \| \triangle\|_F^2,
\end{align*}
where the first equality holds since $\nabla f(T)$ is symmetric,
and the second equality follows from $\big(U\circ U\circ\triangle\big)_{(1)}=\triangle(U\circledast U)^\top$.

For the second class \textbf{(2)}, we have
\begin{align*}
\langle \nabla f(T), -\eta^2(U\circ \triangle\circ \triangle+\cdots)\rangle
=&-3\eta^2\langle \nabla f(T), U\circ \triangle\circ \triangle\rangle
\\
\geq& -3\eta^2\|\nabla f\|\cdot \max_{l\in[r]} \|\u_l\|_2 \cdot \|\triangle\|_F^2\\
=& -\bigg(\underbrace{3\eta\cdot\|\nabla f\|\cdot \max_{l\in[r]} \|\u_l\|_2}_{\leq \frac{1}{L_1} \textbf{(Rule~I)}}\bigg) \bigg(\eta\|\triangle\|_F^2\bigg)
\\
\geq & -\frac{1}{L_1}\cdot\eta\|\triangle\|_F^2,
\end{align*}
where the first inequity follows by using Lemma~\ref{lem:ABC_inner}.

For the third class \textbf{(3)}, we have
\begin{align*}
\eta^3\langle \nabla f(T), \triangle\circ \triangle\circ \triangle\rangle 
\geq & -\eta^3\|\nabla f\| \max_l \|\triangle_l\|_2\cdot \|\triangle\|_F^2
\\
\geq & -\eta^3\|\nabla f\| \bigg(3\|\nabla f\| \max_l\|\u_l\|_2^2\bigg) \|\triangle\|_F^2
\\
=& -\bigg(\underbrace{3\eta\cdot\|\nabla f\|\cdot \max_l\|\u_l\|_2}_{ \leq\frac{1}{L_1}\textbf{(Rule~I)}}\bigg)^2 \bigg(\frac{1}{3}\eta\cdot \|\triangle\|_F^2\bigg)
\\
\geq & -\frac{1}{3L_1^2}\cdot\eta\|\triangle\|_F^2,
\end{align*}
where the first inequality follows from Lemma~\ref{lem:ABC_inner}, and the second inequality follows from Lemma~\ref{lem4:4}.

Combining the above three bounds, we get
$$
\langle \nabla f(T), T-T^+\rangle
\geq \bigg(1-\frac{1}{L_1}-\frac{1}{3L_1^2}\bigg)\cdot \eta\|\triangle\|_F^2.
$$

To bound $M\|T-T^+\|_F^2$, recall that there are three kinds of terms contained in $T-T^+$. Therefore, $\|T-T^+\|_F^2$ has nine terms in total, which can be categorized into two classes. The first class has three terms:
\begin{align*}
&\textbf{(1)}~~~ \|\eta\cdot \triangle\circ U\circ U\|_F^2,
\\
&\textbf{(2)}~~~ \|\eta^2\cdot \triangle\circ \triangle\circ U\|_F^2,
\\
&\textbf{(3)}~~~\|\eta^3\cdot \triangle\circ \triangle\circ \triangle\|_F^2.
\end{align*}
The second class consists of cross terms, i.e., the inner products of any two terms from $\eta\cdot \triangle\circ U\circ U$, $\eta^2\cdot \triangle\circ \triangle\circ U$
and $\eta^3\cdot \triangle\circ \triangle\circ \triangle$. Then, we can bound these cross terms using {\em Cauchy-Schwarz's inequality}. Next, we bound the three terms in the first class in sequence.

First of all, we introduce a lemma that will be frequently used in the remaining part of the proof.
\begin{lem}\cite[Theorem 5.3.4]{horn1994topics}
\label{lem4:odot}
For any semidefinite matrix $A$ and $B$, we have
$$\|A\odot B\|\leq \|A\|\cdot \max_{l} (\diag(B))_{l}. $$
\end{lem}

For the first term \textbf{(1)}, we have
\begin{align*}
\|\eta\cdot \triangle\circ U\circ U\|_F^2
&=\eta^2 \langle U^\top U\odot U^\top U,\triangle^\top \triangle\rangle 
\\
&\leq \eta^2 \|U^\top U\odot U^\top U\|\cdot tr(\triangle^\top \triangle)\\
&\leq \eta^2 \|U^\top U\|\cdot \max_{l\in[r]} (\diag(U^\top U))_{l} \cdot \|\triangle\|_F^2
\\
&=\eta^2  \|U\|^2\cdot \max_{l\in[r]}\|\u_l\|_2^2  \|\triangle\|_F^2\\
&=\bigg(\underbrace{\eta\|U\|^2\cdot \max_{l\in[r]} \|\u_l\|_2^2}_{\leq \frac{1}{ML_2^2} \textbf{(Rule~II)}}\bigg)  \bigg(\eta\|\triangle\|_F^2\bigg)
\\
&\leq \frac{1}{ML^2_2}\cdot\eta\|\triangle\|_F^2,
\end{align*}
where the first inequality follows from Lemma~\ref{lem:aux:von:Neumann}, and the second inequality follows from Lemma~\ref{lem4:odot}.

For the second term \textbf{(2)}, we have
\begin{align*}
\|\eta^2\cdot \triangle\circ \triangle\circ U\|_F^2
=&\eta^4 \langle \triangle^T\triangle\odot U^TU,\triangle^T\triangle\rangle 
\\
\leq& \eta^4 \|\triangle\|_F^2 \cdot \|U\|^2  \bigg(3\|\nabla f\|\cdot \max_l \|\u_l\|_2^2\bigg)^2  \\
=&\bigg(\eta\|\triangle\|_F^2\bigg) \bigg(\underbrace{3\eta\|\nabla f\|\cdot \max_l \|\u_l\|_2}_{\leq\frac{1}{L_1} ~\textbf{(Rule~I)}}\bigg)^2
 \bigg(\underbrace{\eta\|U\|^2\cdot\max_l \|\u_l\|_2^2}_{\leq \frac{1}{ML^2_2} ~\textbf{(Rule~II)}}\bigg)
 \\
\leq& \frac{1}{ML_1^2L^2_2}\cdot\eta\|\triangle\|_F^2,
\end{align*}
where the first inequality follows from Lemmas \ref{lem4:4} and\ref{lem4:odot}, and the fact that $\|U^\top U\|=\|U\|^2$.

For the third term \textbf{(3)}, we have
\begin{align*}
\|\eta^3\cdot \triangle\circ \triangle\circ \triangle\|_F^2
=&\eta^6\langle \triangle^\top \triangle\odot\triangle^\top \triangle,\triangle^\top \triangle \rangle\\
\leq & \eta^6\cdot\|\triangle\|_2^2\cdot\max_l \|\triangle_l\|_2^2\cdot\tmop{tr}(\triangle^\top \triangle)\\
\leq & \eta^6\cdot\big(3\|[\nabla f]_{(1)}\|\cdot \|U\|^2\big)^2\cdot\big(3\|\nabla f\|\max_l\|\u_l\|_2^2\big)^2\cdot \|\triangle\|_F^2\\
=&\bigg(\underbrace{3\eta \|\nabla f\| \max_l \|\u_l\|_2}_{\leq\frac{1}{L_1} ~\textbf{(Rule~I)}}\bigg)^2\bigg(\underbrace{\eta \|U\|^2\cdot \max_l \|\u_l\|_2^2}_{\leq \frac{1}{ML^2_2}  ~\textbf{(Rule~II)}}\bigg)  \bigg(\underbrace{3\|\eta \|[\nabla f]_{(1)}\|\|U\|}_{\leq\frac{1}{L_3} ~\textbf{(Rule~III)}}\bigg)^2 \big(\eta\|\triangle\|_F^2\big)\\
\leq &\frac{1}{ML_1^2L^2_2L_3^2}\cdot\eta\|\triangle\|_F^2,
\end{align*}
where the second line follows from Lemma \ref{lem4:odot}, and the third line follows from Lemmas~\ref{lem4:3} and \ref{lem4:4}.

Recall, we have
\begin{align*}
T-T^+
&= \eta\left(U\circ U\circ\triangle+\cdots\right)
-\eta^2\left(U\circ\triangle\circ\triangle+\cdots\right)
+\eta^3 \triangle\circ \triangle\circ\triangle
\\
&\doteq (A_1+A_2+A_3)-(B_1+B_2+B_3)+C.
\end{align*}
From the bound of the first class terms in $M\|T-T^+\|_F^2$, we can get
\begin{align}\label{lem4:first:class}
\|A_i\|_F^2 &\leq \frac{1}{M}\cdot\frac{1}{L_2^2}\cdot\eta\|\triangle\|_F^2, ~~i=1,2,3,\nonumber\\
\|B_i\|_F^2 &\leq \frac{1}{M}\cdot\frac{1}{L_1^2L^2_2}\cdot\eta\|\triangle\|_F^2 ~~i=1,2,3,\nonumber\\
\|C\|_F^2 &\leq   \frac{1}{M}\cdot\frac{1}{L_1^2L^2_2L_3^2}\cdot\eta\|\triangle\|_F^2.
\end{align}

Now, we are ready to complete the proof by combining the above arguments. In particular, we have
\begin{align*}
\|T-T^+\|_F^2
=&\sum_i \|A_i\|_F^2 + \sum_i \|B_i\|_F^2  +\|C\|_F^2 -\sum_{i,j} \langle A_i,B_j\rangle+\sum_i\langle A_i,C\rangle -\sum_i \langle B_i,C\rangle\\
\leq &\sum_i \|A_i\|_F^2 + \sum_i \|B_i\|_F^2  +\|C\|_F^2+\sum_i\| A_i\|_F\cdot \|C\|_F
+\sum_i \|B_i\|_F\cdot\|C\|_F+\sum_{i,j} \| A_i\|_F\cdot\|B_j\|_F\\
\leq &\left(\frac{3}{L_2^2}+\frac{3}{L_1^2L_2^2}+\frac{1}{L_1^2L_2^2L_3^2}
+\frac{3}{L_1L_2^{2}L_3}+\frac{3}{L_1^2L_2^2L_3}+
\frac{9}{L_1L_2^{2}} \right)\cdot\frac{1}{M}\cdot\eta\|\triangle\|_F^2,
\end{align*}
where the second line follows from the {\em Cauchy-Schwarz's inequality}, and the third inequality follows from~\eqref{lem4:first:class}.
Then, we have
\begin{align}\label{lem4:final:2}
-\frac{M}{2}\|T-T^+\|_F^2
&\geq -\frac{M}{2L_2^2}\frac{1}{M}\bigg(3+\frac{3}{L_1^2}+\frac{1}{L_1^2L_3^2}+\frac{3}{L_1L_3}+\frac{3}{L_1^2L_3}+\frac{9}{L_1} \bigg)\eta\|\triangle\|_F^2\nonumber\\
&= -\frac{1}{2L_2^2}\bigg(3+\frac{3}{L_1^2}+\frac{1}{L_1^2L_3^2}+\frac{3}{L_1L_3}+\frac{3}{L_1^2L_3}+\frac{9}{L_1} \bigg)\eta\|\triangle\|_F^2.
\end{align}

Combining \eqref{lem4:final:1} and \eqref{lem4:final:2}, we get
\begin{align*}
\langle \nabla f(T), T-T^\star \rangle
\geq & \langle \nabla f(T), T-T^+\rangle -\frac{M}{2}\|T-T^+\|_F^2+\frac{m}{2}\|T^\star-T\|_F^2\\
\geq & 
\bigg(
1-\frac{1}{L_1}-\frac{1}{3L_1^2}
-\frac{1}{2L_2^2}
\bigg(3+\frac{3}{L_1^2}+\frac{1}{L_1^2L_3^2}
+\frac{3}{L_1L_3}+\frac{3}{L_1^2L_3}+
\frac{9}{L_1}
\bigg)\bigg)\eta\|\triangle\|_F^2+\frac{m}{2}\|T^\star-T\|_F^2.
\end{align*}
It remains to determine $L_1,L_2,L_3$ precisely such that
\begin{align}\label{lem4:L1:L2:L3}
1-\frac{1}{L_1}-\frac{1}{3L_1^2}
-\frac{1}{2L_2^2}
\bigg(3+\frac{3}{L_1^2}+\frac{1}{L_1^2L_3^2}
+\frac{3}{L_1L_3}+\frac{3}{L_1^2L_3}+
\frac{9}{L_1}
\bigg)\geq\frac{1}{2},
\end{align}
where $L_1,L_2,L_3$ are related with the three stepsize rules.

Next, we simplify  \textbf{Rule~I},  \textbf{Rule~II} and  \textbf{Rule~III} to two new rules  \textbf{Rule~A} and  \textbf{Rule~B}.
Note that  \textbf{Rule~III} implies  \textbf{Rule~I} when $L_3=L_1$. This is due to $\|\nabla f\|\leq\|[\nabla f]_{(1)}\|$.
Moreover, since $ \max_l \|\u_l\|_2 \leq\|U\|$,  \textbf{Rule~II} can be implied by $\eta M\|U\|^4\leq\frac{1}{L_2^2}$.
Thus, by setting $L_3=L_1$, we have the following two new stepsize rules which can imply the original three rules:
\begin{align*}
  \textbf{Rule A}&~~~3L_1\eta \|[\nabla f]_{(1)}\|\|U\|\leq1,\\
  \textbf{Rule B}&~~~\eta L^2_2M\|U\|^4\leq1.
\end{align*}
Note that the requirement \eqref{lem4:L1:L2:L3} now becomes
$$
1-\frac{1}{L_1}-\frac{1}{3L_1^2}
-\frac{1}{2L_2^2}
\bigg(3+\frac{9}{L_1}+\frac{6}{L_1^2}
+\frac{3}{L_1^3}
+\frac{1}{L_1^4}
\bigg)\geq\frac{1}{2}.
$$
Setting $L1=6$ and $L2=3$, we have $\text{LHS}=0.564>\frac{1}{2}.$
Therefore, we set the new stepsize rules as
\begin{align*}
  \textbf{Rule A}&~~~18\eta \|[\nabla f]_{(1)}\|\|U\|\leq1,
 \\
   \textbf{Rule B}&~~~9\eta M\|U\|^4\leq1,
\end{align*}
which can be achieved when the stepsize is chosen as
$
\eta\leq\frac{1}{18\|[\nabla f]_{(1)}\|\|U\|+9M\|U\|^4}.
$

\end{proof}

\end{document}